\documentclass[12pt,letterpaper]{amsart}
\usepackage{euler, epic,eepic,latexsym, amssymb, amscd, amsfonts, xypic, url, color, epsfig}
\input xy
\xyoption{all}
\renewcommand{\familydefault}{ppl}
 \newlength{\baseunit}               % the basic unit length
 \newcount{\numlines}                % depth of picture (in number of lines)
\newcommand{\tpoint}[1]{\vspace{3mm}\par \noindent \refstepcounter{subsection}{\bf \thesubsection.} 
  {\em #1. ---} }
\newcommand{\epoint}[1]{\vspace{3mm}\par \noindent \refstepcounter{subsection}{\thesubsection.} 
  {\em #1.} }

\newcommand{\rK}{{\rm K}}
\newcommand{\rH}{{\rm H}}
\newcommand{\ab}{{\rm ab}}

\newcommand{\Z}{\mathbb{Z}}
\newcommand{\Zhat}{\mathbb{Z}^{\wedge}}

\newcommand{\Q}{\mathbb{Q}}

\newcommand{\kbar}{\overline{k}}

\newcommand{\R}{\mathbb{R}}
\newcommand{\C}{\mathbb{C}}

\newcommand{\proj}{\mathbb P}

\newcommand{\pmR}{\mathbb P^1_{\mathbb{R}} - \{0,1,\infty \}}

\newcommand{\Alb}{\operatorname{Alb}}

\newcommand{\Gal}{\operatorname{Gal}}

\newcommand{\Hom}{\operatorname{Hom}}

\newcommand{\sHom}{\operatorname{\mathcal{H}om}}

\newcommand{\Pic}{\operatorname{Pic}}

\newcommand{\Div}{\operatorname{Div}}

\newcommand{\Spec}{\operatorname{Spec}}
\newcommand{\Sheafspec}{\underline{\operatorname{Spec}}}
\newcommand{\Sheafproj}{\underline{\operatorname{Proj}}}

\newcommand{\Sym}{\operatorname{Sym}}

\newcommand{\Ker}{\operatorname{Ker}}

\newcommand{\Image}{\operatorname{Image}}

\newcommand{\topp}{\textrm{top}}

\newcommand{\hidden}[1]{\footnote{Hidden:  #1}}
\renewcommand{\hidden}[1]{}

\begin{document}
\pagestyle{plain}
\title{$2$-Nilpotent Real Section Conjecture}

\author{Kirsten Wickelgren}\thanks{Supported by an American Institute of Mathematics $5$-year fellowship.}
\address{Dept. of Mathematics, Harvard University, Cambridge~MA}
%\email{vakil@math.stanford.edu}
\date{Tuesday, January 31, 2012.}
\begin{abstract}
We show a $2$-nilpotent section conjecture over $\R$: for a geometrically connected curve $X$ over $\R$ such that each irreducible component of its normalization has $\R$-points, $\pi_0(X(\R))$ is determined by the maximal $2$-nilpotent quotient of the fundamental group with its Galois action, as the kernel of an obstruction of Jordan Ellenberg. This implies that for $X$ smooth and proper, $X(\R)^{\pm}$ is determined by the maximal $2$-nilpotent quotient of $\Gal(\C(X))$ with its $\Gal(\R)$ action, where $X(\R)^{\pm}$ denotes the set of real points equipped with a real tangent direction, showing a $2$-nilpotent birational real section conjecture. 
\end{abstract}
\maketitle
%\tableofcontents

{\parskip=12pt % closing bracket is just before the bibliography 

\section{Introduction}

Grothendieck's section conjecture predicts that the rational points of hyperbolic curves over finitely generated fields are determined by their \'etale fundamental groups. Let $X$ denote a geometrically connected, finite type scheme over a characteristic $0$ field $k$, equipped with a geometric point $b: \Spec \Omega \to X$\hidden{, i.e. $\Omega$ is an algebraically closed field}. Let $\kbar$ denote the algebraic closure of $k$ in $\Omega$. There is a canonical lift of $b$ to $X_{\kbar} = X \times_{\Spec k} \Spec \kbar$ and an exact sequence of \'etale fundamental groups \begin{equation}\label{htpyes}1 \to \pi_1(X_{\kbar}, b) \to \pi_1(X, b) \to G_k \to 1,\end{equation} where $G_k = \Gal (\kbar/ k)$ is the absolute Galois group of $k$, and all fundamental groups are based at the geometric points naturally associated to $b$ \cite[IX Thm 6.1]{sga1}. A rational point $x: \Spec k \to X$ induces a map $\pi_1(x): G_k \to \pi_1(X, x)$, where $\pi_1(X,x)$ denotes the \'etale fundamental group of $X$ based at the geometric point $\Spec \kbar \to \Spec k \to X$ associated to $x$. View $x$ and $b$ as geometric points of $X_{\kbar}$ and choose a path between them, where path means a natural transformation between the associated fiber functors, giving a path between $x$ and $b$ in $X$ and an isomorphism $\pi_1(X,x) \cong \pi_1(X,b)$ respecting the projections to $G_k$. Composing $\pi_1(x)$ with this isomorphism $\pi_1(X,x) \cong \pi_1(X,b)$ produces a section $s: G_k \to \pi_1(X,b)$ of \eqref{htpyes}, and a different choice of path will change the section  to $g \mapsto \lambda s(g) \lambda^{-1}$ for some $\lambda$ in $\pi_1(X_{\kbar}, b)$. Sections obtained from $s$ in this way are said to be conjugate. Let $\mathscr{S}_{\pi_1(X/k)}$ denote the conjugacy classes of sections of \eqref{htpyes}, and let $X(k)$ denote the set of $k$-points of $X$. Given $X$ and $b$ as above, let $\kappa$ be the map $$\kappa: X(k) \to \mathscr{S}_{\pi_1(X/k)}$$ just constructed.  For $X$ a smooth, proper curve of genus $> 1$ over a finitely generated field, Grothendieck's section conjecture, which is unknown, is that $\kappa$ is a bijection. 

For $k=\R$, the map $\kappa$ factors through $\pi_0(X(\R))$, and the real section conjecture, saying that $$\kappa:  \pi_0(X(\R)) \to \mathscr{S}_{\pi_1(X/\R)}$$ is a bijection is proven, but non-trivial \cite{Mochizuki_real_SC} \cite{Sullivan} \cite{Miller_Sullivan} \cite{Carlsson_Sullivans_Conj} \cite{Pal}. This paper proves a $2$-nilpotent real section conjecture, determining $\pi_0(X(\R))$ from the maximal $2$-nilpotent quotient of $\pi_1(X_{\C})$ with its $G_{\R}$-action.

For a (profinite) group $\pi$, let $\pi = [\pi]_1> [\pi]_2 > [\pi]_3 > \ldots$ denote the lower central series of $\pi$, i.e. $[\pi]_{n+1} = [[\pi]_n,\pi]$ (respectively $[\pi]_{n+1}= \overline{[[\pi]_n,\pi]}$) is (the closure of) the subgroup generated by commutators of elements of $[\pi]_n$ and $\pi$. Pushing out \eqref{htpyes} by the quotient $\pi_1(X_{\kbar}, b) \to \pi_1(X_{\kbar}, b)/[\pi_1(X_{\kbar}, b)]_n$ yields an exact sequence \begin{equation}\label{htpyesn=2} 1 \to \pi_1(X_{\kbar}, b)/[\pi_1(X_{\kbar}, b)]_n \to \pi_1(X, b)/[\pi_1(X_{\kbar}, b)]_n \to G_k \to 1.\end{equation} Let $\kappa^{\ab}$ be the map taking a $k$-point $x$ of $X$ to the section of \eqref{htpyesn=2} for $n=2$ determined by $\kappa(x)$.

Define {\em curve} to mean a pure dimension $1$, finite type scheme over a field. A curve $X$ over $k$ will be said to be {\em based} if it is equipped with a choice of a geometric point whose image is a $k$-point of $X$ or which is associated to a $k$-tangent vector based at a smooth point of a compactification of $X$ as described in \cite[\S 15]{Deligne} \cite{Nakamura} \cite[12.2.1]{PIA}.\hidden{ The \'etale fundamental group determines a functor from based curves over $k$ to profinite groups over $G_k$.} The complex analytic space $X(\C)$ associated to a based curve $X$ over $\R$ has a distinguished point or tangent vector based at a smooth point of a compactification\hidden{, which is fixed by the associated $G_{\R}$-action}, allowing us to apply the topological or orbifold fundamental group functors to $X(\C)$ or $X(\C)/ G_{\R}$, giving maps $\kappa$ and $\kappa^{\ab}$ as above. 

In the following theorem, $X$ is a based curve over $\R$, $\pi$ denotes either the \'etale or topological fundamental group of $X_{\C}$ or $X(\C)$ respectively, and $\pi_1(X)$ denotes either the \'etale or orbifold fundamental group of $X$ or $X(\C)/ G_{\R}$ respectively.

\tpoint{Theorem}\label{2nilsectionconjecture}{\em Let $X$ be a geometrically connected, based curve over $\R$, such that each irreducible component of its normalization has $\R$-points. Then $\kappa^{\ab}$ is a natural bijection from $\pi_0(X(\R))$ to conjugacy classes of sections of $$ 1 \to \pi/[\pi]_2 \to \pi_1(X)/[\pi]_2 \to G_{\R} \to 1 $$ which lift to sections of $$1 \to \pi/[\pi]_3 \to \pi_1(X)/[\pi]_3 \to G_{\R} \to 1.$$}

Note that the assumption that $X$ is based gives \eqref{htpyes} a splitting, and that this implies that Theorem \ref{2nilsectionconjecture} says that the $2$-nilpotent quotient $\pi/[\pi]_3$ of $\pi$ with its $G_{\R}$-action determines the connected components of $X(\R)$. The real section conjecture shows that $\pi$ with its $G_{\R}$-action determines the connected components of $X(\R)$ when the topological space $X(\C)$ is a $\rK(\pi,1)$, which is the case precisely when no component of the normalization of $X_{\C}$ is $\proj^1$ -- see Remark \ref{relation_condition_normalization_components_not_pointless_to_condition_X(C)_is_K(pi,1)}. The proof of Theorem \ref{2nilsectionconjecture} given below is independent of the real section conjecture, although assuming it, one would be saved the trouble of proving Proposition \ref{kappaJaciso}.

For $X$ smooth and proper, Theorem \ref{2nilsectionconjecture}  applied to smaller and smaller Zariski opens of $X$ shows that $X(\R)$ is determined by the maximal $2$-nilpotent quotient of the absolute Galois group of the function field $\C(X)$ of $X_{\C}$ with its $G_{\R}$-action. \hidden{More specifically, for a geometrically integral curve $X$ over $k$, and a field extension $L/k$, let $L(X)$ denote the function field of $X_{L}$. For a smooth, proper, connected, based curve $X$ over $\R$,\hidden{claim: connected and based implies geometrically connected pf: let $f: X \to \Spec k$ denote the structure map of a connected, based curve over $k$. We have $\overline{k} \otimes f_* \mathcal{O}_X= (f \otimes \overline{k})_* \mathcal{O}_{X \otimes \overline{k}}$ as can be seen directly from the \v{C}ech complexes computing both sides. Since $X$ is proper, $f_* \mathcal{O}_X$ is a finite $k$-module. Since $X$ is connected and smooth, $X$ is irreducible, whence $f_* \mathcal{O}_X = \rH^0(X, \mathcal{O}_X)$ is an integral domain. It follows that $f_* \mathcal{O}_X$ is a field. Since $X$ is based, we must have $f_* \mathcal{O}_X = k$, whence the global sections of $\mathcal{O}_{X \otimes \overline{k}}$ is $\overline{k}$, whence $X \otimes \overline{k}$ is connected.}} Let $X(\R)^{\pm}$ denote the set of real points of $X$ equipped with a real tangent direction, i.e. a vector in the tangent space of the smooth $1$-manifold $X(\R)$ up to scaling by elements of $\R_{>0}$. The notation $X(\R)^{\pm}$ is meant to indicate that after orienting $X(\R)$, the two tangent directions associated to each element of $X(\R)$ consist of the direction distinguished by the orientation and its negative. For any Zariski open $U$ of $X$, there is a map $X(\R)^{\pm} \to \pi_0(U(\R))$ given by taking a tangent direction to the connected component it is pointing towards. Note that the resulting map $X(\R)^{\pm} \to \varprojlim_{U} \pi_0(U(\R))$ is a bijection. It follows that a corollary of the $2$-nilpotent real section conjecture is that $X(\R)^{\pm}$ is determined by $G_{\R(X)}/[G_{\C(X)}]_3 \to G_{\R}$ (\S 4 Corollary \ref{corbir2nilSC}).

\tpoint{Corollary}\label{corbir2nilSCintro}{\em Let $X$ be a smooth, proper, connected curve over $\R$ equipped with a chosen element of $X(\R)^{\pm} \neq \emptyset$. There is a natural bijection between $X(\R)^{\pm}$ and the conjugacy classes of sections of $$ 1 \to G_{\C(X)}/[G_{\C(X)}]_2 \to G_{\R(X)}/[G_{\C(X)}]_2 \to G_{\R} \to 1 $$ which lift to sections of $$1 \to G_{\C(X)}/[G_{\C(X)}]_3 \to G_{\R(X)}/[G_{\C(X)}]_3 \to G_{\R} \to 1.$$ }

The real section conjecture and its $2$-nilpotent version are closely related to Sullivan's conjecture, as we now discuss, first introducing some notation. This also helps summarize the proof of Theorem \ref{2nilsectionconjecture}, which we do below.

Let $G= \Z/2$ and $EG$ denote a contractible topological space with a free action of $G$. For a sufficiently well-behaved topological space $Y$ with a $G$-action, e.g. $Y$ a $G$-CW complex, let $\operatorname{Map}(E G , Y)$ denote the function space of continuous maps $E G \to Y$ equipped with the $G$ action given by $g f = gfg^{-1}$. The homotopy fixed points of $G$ on $Y$ are defined $Y^{h G} = \operatorname{Map}(E G, Y)^{G}$ and there is a canonical map $Y^G \to Y^{hG}$ from the fixed points to the homotopy fixed points induced by the $G$-equivariant map from $E G$ to the point.  

Let $\mathscr{S}_{\pi_1(Y/G)}$ denote the conjugacy classes of sections of \begin{equation}\label{orbifold_pi_1_sequence} 1 \to \pi \to \pi_1(Y) \to G \to 1 \end{equation} where $\pi$ denotes the topological fundamental group of $Y$, based at some point not included in the notation, and $\pi_1(Y) $ denotes the orbifold fundamental group, which can be identified with the topological fundamental group of $EG \times_G Y$ or with the group of automorphisms of the universal cover of $Y$ lying over an automorphism of $Y$ induced by an element of $G$. There is a natural map $\pi_0(Y^{hG}) \to \mathscr{S}_{\pi_1(Y/G)}$, which is a bijection if $Y$ is a $\rK(\pi,1)$.

For $X$ a geometrically connected, finite type scheme over $\R$, the map $\kappa$ for the \'etale fundamental group is the composition \begin{equation}\label{kappa_as_composition} \pi_0 (X(\C)^{G_{\R}}) \to  \pi_0(X(\C)^{h G_{\R}}) \to \mathscr{S}_{\pi_1(X(\C)/G_{\R})} \to \mathscr{S}_{\pi_1(X/\R)}\end{equation} where the last map is induced by the canonical isomorphism from the profinite completion of \eqref{orbifold_pi_1_sequence} to \eqref{htpyes} \cite[XII Cor 5.2]{sga1}, and the map $\kappa$ for the topological fundamental group is the composition of the first two maps of \eqref{kappa_as_composition}. For $X(\C)$ a $\rK(\pi,1)$, as in the section conjecture, the second map is a bijection. 

The Sullivan conjecture \cite{Sullivan}, proven by Miller \cite{Miller_Sullivan}, Dwyer-Miller-Neisendorfer \cite{DMR_Sullivan}, Carlsson \cite{Carlsson_Sullivans_Conj}, and Lannes \cite{Lannes_FBVX}, shows that the first map of \eqref{kappa_as_composition} is a bijection. Precisely, it says that the natural map from the $p$-completion of the fixed points to the homotopy fixed points of the $p$-completion is a weak equivalence for a finite $p$-group $G$ acting on a finite $G$-CW complex, but proven at the same time is the fact that applying $\pi_0$ to $Y^{G} \to Y^{h G}$, as in the first map, is a bijection \cite[Theorem B (a)]{Carlsson_Sullivans_Conj}. So if one overlooks the map $\mathscr{S}_{\pi_1(X(\C)/G_{\R})} \to \mathscr{S}_{\pi_1(X/\R)}$ comparing the topological to the \'etale fundamental group, i.e. if one uses $\kappa$ for the topological fundamental group, the real section conjecture is $\pi_0$ of Sullivan's conjecture applied to a $\rK(\pi,1)$. Also see \cite{Pal} for a nice proof of the real section conjecture in the \'etale and topological case which does not appeal to Sullivan's conjecture.

The proof of Theorem \ref{2nilsectionconjecture} can be summarized as follows. Let $\pi_0 (X(\R))^-$ denote $\pi_0 (X(\R))$ with the component containing the base point removed. Let $\pi_0((-)^{G_{\R}})$ denote the functor taking $G_{\R}$-fixed points and then applying $\pi_0$. Applying $\pi_0((-)^{G_{\R}})$ to the map  $X(\C) \to \Sym^{\infty} X(\C)$ from $X(\C)$ to its infinite symmetric product expresses $\pi_0 (X(\R))^-$ as a basis for the $\Z/2$-vector space $\pi_0( (\Sym^{\infty} X(\C))^{G_{\R}})$ (Proposition \ref{pi_0SymXG_freeVS}). For a smooth curve $X$ with generalzed Jacobian $\Pic^0 X^+$, there is a natural map $\Sym^{\infty} X(\C) \to \Pic^0 X^+ (\C)$ which is an affine or projectivized vector bundle, whence a bijection after applying $\pi_0((-)^{G_{\R}})$ (Proposition \ref{pi0(Syminfty_to_Pic)GR_iso}). The conjugacy classes of sections of \eqref{htpyesn=2} for $n=2$ are canonically identified with the connected components of the homotopy fixed points $(\Pic^0 X^+ (\C))^{h G_{\R}}$ as above, since $\Pic^0 X_{\C}^+$ is a $\rK(\pi/[\pi]_2,1)$, but $\kappa$ also identifies them with the connected components of the fixed points of  $\Pic^0 X^+ (\C)$, as is shown without appealing to Sullivan's conjecture (Proposition \ref{kappaJaciso}). In total, it follows that $\pi_0 (X(\R))^-$ is a $\Z/2$-basis for the conjugacy classes of sections of \eqref{htpyesn=2} for $n=2$, and in fact this holds without the assumption that $X$ is smooth (Proposition \ref{kappaab_free_vs}).  

By a standard interpretation of group cohomology, the conjugacy classes of sections of \eqref{htpyesn=2} for $n=2$ are identified with $\rH^1(G_{\R}, \pi/[\pi]_2)$. The obstruction to lifting to a section of \eqref{htpyesn=2} for $n=3$ is a map $$\delta_2: \rH^1(G_{\R}, \pi/[\pi]_2) \to \rH^2(G_{\R}, [\pi]_2/[\pi]_3),$$ which is quadratic with associated bilinear form $$\rH^1(G_{\R}, \pi/[\pi]_2) \wedge \rH^1(G_{\R}, \pi/[\pi]_2) \to \rH^2(G_{\R}, [\pi]_2/[\pi]_3 ) $$ induced by the cup-product and the commutator pairing $$\pi/[\pi]_2 \wedge \pi/[\pi]_2 \to [\pi]_2/[\pi]_3,$$ as follows from a result of Zarkhin \cite[p 242]{Zarkhin}. We show the associated bilinear form is injective (Lemmas \ref{cupwedgeinjectiveZ/2} and \ref{[-,-]*injective}), which implies Theorem \ref{2nilsectionconjecture} (Theorem \ref{kerdelta2=Imkappaab}).

It is not true in general that Sullivan's conjecture holds for the infinite symmetric product of a finite $G$-CW complex, nor that \begin{equation}\label{pi0SymintyGtohG} \pi_0 ((\Sym^{\infty} Y)^{G}) \to \pi_0 ((\Sym^{\infty} Y)^{h G}) \end{equation} is a bijection. An interesting consequence of the proof of the  $2$-nilpotent real section is that for $Y = X(\C)$ with $X$ a real based algebraic curve such that each irreducible component of its normalization has $\R$-points, the map \eqref{pi0SymintyGtohG} is a bijection. See Remark \ref{pi0(SymGtoSymhG)bij}. 

$\delta_2$ was studied by Jordan Ellenberg as an obstruction to rational points of a curve's Jacobian lying in the image of an Abel-Jacobi map \cite{Ellenberg_2_nil_quot_pi}, and also studied by Zarhin \cite{Zarkhin}. Theorem \ref{2nilsectionconjecture} was guessed by Jordan Ellenberg in the proper, smooth case, as he told me after I had observed it held in several examples, and it can naturally be expressed in terms of his ideas in \cite{Ellenberg_2_nil_quot_pi}: a smooth based curve embeds into its generalized Jacobian by its Abel-Jacobi map. Those rational points $y$ of the Jacobian which are the image of a point of the curve satisfy the condition that $\kappa(y)$ lifts to a section of \eqref{htpyes} where $X$ denotes the curve. Filtering $\pi_1(X_{\C})$ by its lower central series provides a series of obstructions, the first of which is the quadratic form, here denoted $\delta_2$ and discussed in \S \ref{2nilobstruction_section}.   

We give a topological interpretation of Ellenberg's point of view in Section \ref{top2nilapprox}, constructing a diagram of finite $G_{\R}$-CW complexes \begin{equation}\xymatrix{ & \Alb_2 \ar[d]^q \\
X(\C) \ar[ur] \ar[r] & \Alb_1} \end{equation} for an arbitrary geometrically connected curve $X$ over $\R$ with a chosen real base point, such that $\Alb_2$ is a $\rK(\pi_1(X(\C))/[\pi_1(X(\C))]_3, 1)$, $\Alb_1$ is a $\rK(\pi_1(X(\C))/[\pi_1(X(\C))]_2, 1)$, $q$ is a fiber bundle, and all maps induce the obvious quotient maps on topological fundamental groups. Sullivan's conjecture gives an equivalence between Theorem \ref{2nilsectionconjecture} and the statement that the connected components of the real points of the curve are those of the abelian approximation $\Alb_1$ which can be lifted to the $2$-nilpotent approximation $\Alb_2$. See Theorem \ref{Top_2nilSC}.

{ \bf Relation to other work:} Grothendieck's section conjecture is part of his anabelian conjectures predicting that certain schemes are determined by their \'etale fundamental groups. Birational variants of the anabelian conjectures replace $\pi_1^{et}$ by the absolute Galois group of the function field. There has been considerable work describing varieties using small quotients of their fundamental groups or the Galois groups of their function fields. Pop has shown a meta-abelian birational section conjecture over $p$-adic fields \cite{Pop_bir_p-adic_sC}. Bogomolov and Tschinkel developed an approach to recognize the function field of certain varieties using the $2$-nilpotent quotient of the absolute Galois group. Work of Bogomolov, Pop, and Tschinkel shows that it is possible to recover the function field of certain varieties of dimension $\geq 2$ over algebraically closed fields from the $2$-nilpotent quotient of the absolute Galois group \cite{Bogomolov_On_two_conj_bir_anab} \cite{Bogomolov_Ab_sub_Gal_grps} \cite{MR1977585} \cite{BT_Recon_fun_fields_geom} \cite{Pop28} \cite{Pop29} -- see \cite{BTSurvey} and \cite{LangMemorial} for more discussion. There is also interesting work limiting when such minimalistic anabelian results can hold. Yuichiro Hoshi has found examples where any section of a pro-$p$ homotopy exact sequence of the Jacobian lifts to a section of a pro-$p$ homotopy exact sequence of the curve \cite{Hoshi_non-geom_pro-pGalois_sec}\hidden{In preprint, this was Theorem 3.5, Corollary 3.6}.

{\bf Acknowledgments: } I wish to thank Jordan Ellenberg for sharing \cite{Ellenberg_2_nil_quot_pi} with me and for suggesting that \cite{GrossHarris} could be used to prove the main result of this paper in the smooth, proper case. I also wish to thank Gunnar Carlsson for recognizing the connection between this problem and Sullivan's conjecture. It is a pleasure to express my deepest gratitude and admiration to both Carlsson and Ellenberg. I thank an anonymous reviewer of a previous version of this paper for Lemma \ref{cupwedgeinjectiveZ/2} and other helpful suggestions. I would also like to thank Florian Pop for suggesting that Theorem \ref{2nilsectionconjecture} could be used to show a birational result, as well as for many helpful comments.  

\section{Abelian approximation}

We compute $\kappa^{\ab}$, which is considered as an abelian approximation to $\pi_0(X(\R))$.

Let $X$ be a geometrically connected scheme with geometric point $b$ and \'etale fundamental group $\pi$. The base point $b$ determines a natural bijection $\mathscr{S}_{\pi_1(X/k)} = \rH^1(G_k, \pi)$. Under this identification, $\kappa(x)$ is represented by the cocycle $$G_k \ni g \mapsto \gamma^{-1} (g \gamma) \in \pi,$$ where $\gamma$ is a path from the base point to $x \in X(k)$, and composition in the fundamental group is written so that $\gamma^{-1} (g \gamma)$ is the loop starting at the base point obtained by first traversing $g \gamma$ and then traversing $\gamma^{-1}$.  The analogous statements of course hold for the topological space $X(\C)$.

Let $\mathcal{I}$ denote the forgetful functor from vector spaces over $\Z/2$ to pointed sets, sending a vector space to its underlying set, pointed by the identity. Let $\mathcal{V}$ denote its left adjoint, called the {\em free vector space} on the pointed set. The {\em unit} is the canonical map of pointed sets $( S, s_0 ) \rightarrow \mathcal{I} \mathcal{V} ( S, s_0 )$, so $\mathcal{V}(S, s_0)$ has basis $S -\{s_0\}$, and the unit map sends $S - \{s_0\}$ to this basis and $s_0$ to $0$.

Note that for a based curve $X$ over $\R$, the set $\pi_0(X(\R))$ is naturally pointed, as the base point's image either lies in a particular connected component or the associated tangent vector points towards one. 

\tpoint{Proposition}\label{kappaab_free_vs}{\em Let $X$ be a geometrically connected, based curve over $\R$, such that each irreducible component of its normalization has $\R$-points. Then $$\kappa^{\ab}: \pi_0(X(\R)) \rightarrow \rH^1(G_{\R}, \pi^{\ab})$$ is canonically isomorphic to the unit of the adjunction $(\mathcal{V}, \mathcal{I})$ on the pointed set $\pi_0(X(\R))$.}

In Proposition \ref{kappaab_free_vs}, $\pi$ can denote either the \'etale or topological fundamental group of $X_{\C}$ or $X(\C)$, respectively.

The numerical version of Proposition \ref{kappaab_free_vs} for $X$ smooth and proper saying that $$2^{\vert \pi_0(X(\R))\vert - 1} = \vert \rH^1(G_{\R}, \pi^{\ab})\vert$$ follows from combining \cite[Prop 1.3, Prop 3.2]{GrossHarris} with the Kummer exact sequence of $\Pic^ 0 X$. There is a similar numerical computation of $\vert \pi_0(X(\R)) \vert$ in terms of $\Z/2$-homology of $X$ in \cite[Prop 4.4]{GrossHarris}. The definition of $\kappa^{\ab}$ allows for the natural version above. 

\begin{proof}
We first reduce to the case where $X$ is smooth. Let $f: \tilde{X} \to X$ be the normalization of $X$. Note that $\tilde{X}$ is a disjoint union of smooth curves satisfying the hypotheses of the proposition. There is a finite $G_{\R}$-equivariant set $D \subset X(\C)$ such that $X(\C)$ is homeomorphic to the push-out $ \tilde{X}(\C) \coprod_{\tilde{D}} D $ where $\tilde{D} = f^{-1} (D)$ and the push-out is taken with respect to $f: \tilde{D} \to D$ and the inclusion $\tilde{D} \subset \tilde{X}(\C)$. It follows that \begin{equation}\label{pi0XXtildecoprod}\pi_0 X(\R) \cong  \pi_0 \tilde{X}(\R) \coprod_{\tilde{D}^{G_{\R}}} {D}^{G_{\R}}.\end{equation}

 By Mayer-Vietoris, there is an exact sequence \begin{equation}\label{H1ES} 0 \to \rH_1(\tilde{X}(\C)) \to \rH_1(X(\C)) \to \oplus_{\tilde{D}} \Z \to \oplus_{D\coprod \pi_0 (\tilde{X}(\C)) } \Z \to \Z \to 0 \end{equation} of $\Z[G_{\R}]$-modules, where $ \rH_1(-)$ denotes singular homology with $\Z$-coefficients. For the \'etale case, substitute $\Zhat$ for $\Z$ and the abelianization of $\pi_1$ for $\rH_1(-)$. The sequence \eqref{H1ES} remains exact and the following argument is valid with these substitutions.

Let $\tau$ denote complex conjugation. Consider the double complex $$(E_{ij}, d_{ij}: E_{ij} \to E_{i+1, j}, \delta_{ij}: E_{ij} \to E_{i, j + 1}  )$$ with identical rows $(E_{* j}, d_{*j})$ equal to \eqref{H1ES} and differentials $\delta_{i, 2j} = 1 - \tau$ and $\delta_{i, 2j} = 1 + \tau$. This double complex gives a spectral sequence $(E^r_{ij}, D^r_{ij}: E^r_{ij} \to E^r_{i+r, j - r + 1} )$, $r=0,1,\ldots$ and $d^0_{ij} = \delta_{ij}$. This spectral sequence converges with $E^5_{ij} = E^{\infty}_{ij} = 0$, because \eqref{H1ES} is exact.

Since $\oplus_{\tilde{D}} \Z$ and $\oplus_{D\coprod \pi_0 (\tilde{X}(\C)) } \Z$ are of the form $\Z[G_{\R}]^a \oplus \Z^b$, it follows that $(E^1_{*j}, d^1_{*j})$ for $j = 2n,2n+1$ is $$\xymatrix { \rH^1(G_{\R}, \rH_1(\tilde{X}(\C))) \ar[r] & \rH^1(G_{\R}, \rH_1(X(\C))) \ar[r] & 0 \ar[r] &0 \ar[r] & 0 \\
 \hat{\rH}^0(G_{\R}, \rH_1(\tilde{X}(\C))) \ar[r] & \hat{\rH}^0(G_{\R}, \rH_1(X(\C))) \ar[r] & \oplus_{\tilde{D}^{G_{\R}}} \Z/2 \ar[r] & \oplus_{D^{G_{\R}}\coprod \pi_0 (\tilde{X}(\C))} \Z/2 \ar[r] & \Z/2  } $$ where $\hat{\rH}^0 = \rH^2$ is Tate cohomology. The map $\kappa^{\ab}$ applied to $D^{G_{\R}}$ and the images of chosen base points for each component of $\tilde{X}(\C)$ gives a splitting of $$d^2_{11}:   \rH^1(G_{\R}, \rH_1(X(\C))) \to E^2_{03},$$ inducing an exact sequence \begin{equation*}\label{ESH1GH1} 0 \to \oplus_{\tilde{D}^{G_{\R}}} \Z/2  \to \rH^1(G_{\R}, \rH_1(\tilde{X}(\C))) \oplus ( \oplus_{D^{G_{\R}}\coprod \pi_0 (\tilde{X}(\C))} \Z/2) \to \rH^1(G_{\R}, \rH_1(X(\C))) \oplus \Z/2 \to 0 \end{equation*} which is compatible with $\kappa^{\ab}$ applied to the coproduct decomposition \eqref{pi0XXtildecoprod} of $\pi_0(X(\R))$ and the resulting short exact sequence $$0 \to \oplus_{\tilde{D}^{G_{\R}}} \Z/2  \to \oplus_{ \pi_0(\tilde{X}(\R))  \coprod D^{G_{\R}} }\Z/2  \to \oplus_{\pi_0(X(\R))} \Z/2 \to 0.$$ It follows that it is sufficient to prove the claim for each connected component of $\tilde{X}$, i.e. we may assume that $X$ is smooth.

We may also assume that the base point $b$ of $X$ is a geometric point whose image is an $\R$-point, i.e $b$ is not tangential: for $b$ and $b'$ different choices of base point, $$\kappa_b^{\ab} (x)= \kappa_{b'}^{\ab} (x) +  \kappa_b^{\ab}(b')$$ for all $x$ in $\pi_0(X(\R))$. In particular, $\kappa^{\ab}_b$ is canonically isomorphic to $\kappa^{\ab}_{b'}$ for any $b'$ which determines the same path component of $X(\R)$.

Since $X$ is smooth, $X$ embeds into its generalized Jacobian: for $X$ smooth, non-proper, let $X^+$ denote the coproduct in schemes $X^c \coprod_{X^c - X} \Spec \R$ where $X^c$ denotes the smooth compactification of $X$, and let $X^+ = X$ for $X$ smooth, proper. In other words, $X^+$ is the one-point compactification formed by crushing $X^c -X$ to a point. Sending a point $x$ of $X$ to the invertible sheaf of rational functions on $X^+$ with at worst a simple pole at $x$ determines a map $X \to \Pic^1 X^+$ from $X$ to the degree $1$ Picard scheme of $X^+$. Translation by the $\R$-point of $\Pic^1 X^+$ equal to the image of $b$ gives an isomorphism $\Pic^n X^+ \cong \Pic^0 X^+$ for all $n$. The resulting map $$\alpha: X \to \Pic^0 X^+$$ is the Abel-Jacobi embedding of $X$ into its generalized Jacobian.

The Abel-Jacobi map induces $\pi_1(\alpha): \pi \to \pi_1( \Pic^0 X^+_{\C})$ which is the abelianization in either the \'etale or topological case \cite[Prop A.8 (iii)]{Mochizuki_Top_Ab_An_Geom_I}. By functoriality of $\kappa$, the diagram $$\xymatrix{\pi_0(X(\R))\ar[d]^{\kappa} \ar[rr]^{\pi_0(\alpha(\R))}&& \pi_0(\Pic^0 X^+ (\R)) \ar[d]^{\kappa_J }\\
 \rH^1(G_{\R}, \pi ) \ar[rr]^{\pi_1(\alpha)_*} && \rH^1(G_{\R}, \pi_1( \Pic^0 X^+_{\C}) )}$$ commutes, where $\kappa_J$ is based at the identity of $\Pic^0 X^+$, giving an isomorphism between $\kappa^{\ab}$ and the composition $$\xymatrix{\pi_0(X(\R)) \ar[rr]^{\pi_0(\alpha(\R))} &&  \pi_0(\Pic^0 X^+ (\R)) \ar[r]^{\kappa_J } & \rH^1(G_{\R}, \pi_1( \Pic^0 X^+_{\C}) )}.$$

Proposition \ref{kappaab_free_vs} follows from the fact that $\kappa_J$ is an isomorphism (Proposition \ref{kappaJaciso}) and $\pi_0(\alpha(\R))$ is canonically isomorphic to the unit of the adjunction $(\mathcal{V}, \mathcal{I})$ on the pointed set $\pi_0(X(\R))$, as follows from Proposition \ref{pi_0SymXG_freeVS} and \ref{pi0(Syminfty_to_Pic)GR_iso}.

\end{proof}

Let $X$ be a smooth geometrically connected, based curve over $\R$, and choose a geometric point over the identity of $\Pic^0 X^+_{\C}$, giving $\kappa_J: \pi_0(\Pic^0 X^+ (\R)) \to \rH^1(G_{\R}, \pi_J )$, where $\pi_J$  denotes either the \'etale or topological fundamental group of $\Pic^0 X^+_{\C}$.  

 \tpoint{Proposition}\label{kappaJaciso}{\em The map $\kappa_J$ is an isomorphism of $\Z/2$-vector spaces.}
 
 The content of the proof of Proposition \ref{kappaJaciso} is identical in the \'etale and topological setting using \cite[\S 3.3, 4.1]{VWfund_grp}. We give the \'etale proof.
 
 \begin{proof}
 Let $\tilde{J}$ denote the universal cover of $J = \Pic^0 X^+$, which is automatically an abelian group. The canonical exact sequence $$0 \to \pi_J \to \tilde{J} \to J_{\C} \to 0 $$ gives a short exact sequence of abelian groups with $G_{\R}$-action $$0 \to \pi_J \to \tilde{J}(\C) \to J(\C) \to 0 .$$ Applying Tate cohomology gives the exact sequence \begin{equation}\label{TatecohHGLESpiJtildeJJ} \hat{\rH}^0(G_{\R}, \tilde{J}(\C)) \rightarrow \hat{\rH}^0(G_{\R}, J(\C)) \rightarrow  \hat{\rH}^1(G_{\R}, \pi_J ) \rightarrow \hat{\rH}^1(G_{\R}, \tilde{J}(\C)).\end{equation} 
 
 Since $\tilde{J} \to J \stackrel{[2]}{\to} J$ and $\tilde{J} \to J$ are simply connected covering spaces of $J$, where $[2]: J \to J$ denotes multiplication by $2$, they are isomorphic  \cite[Prop. 3.1, Thm 3.1]{VWfund_grp}, whence multiplication by $2$ is an isomorphism $[2]: \tilde{J} \to \tilde{J}$ and $\hat{\rH}^i(G_{\R}, \tilde{J})=0$ for all $i$.
 
It is straight-forward to verify that the diagram $$\xymatrix{\hat{\rH}^0(G_{\R}, J(\C)) \ar[r]^{\cong} &  \hat{\rH}^1(G_{\R}, \pi_J ) \\ J(\R) \ar[u] \ar[r] & \pi_0(J(\R)) \ar[u]_{\kappa_J} }$$ is commutative. Since $J(\R) \to \hat{\rH}^0(G_{\R}, J(\C))$ is surjective, $\kappa_J$ is surjective. By \cite[Prop 1.3]{GrossHarris}, $\hat{\rH}^0(G_{\R}, J(\C))$ and $\pi_0(J(\R))$ have the same cardinality, so $\kappa_J$ is bijective.
 \end{proof} 
 
\epoint{Remark} Proposition \ref{kappaJaciso} also follows from the real section conjecture. %The above proof is included in order to give a proof of the $2$-nilpotent real section conjecture which is independent of the real section conjecture or Sullivan's conjecture.

The unit of the adjunction $(\mathcal{V}, \mathcal{I})$ on $\pi_0(X(\R))$ is computed by the following proposition, whose proof is is omitted here, but is essentially contained in the proof of Proposition 3.2 in \cite{GrossHarris}, which Gross and Harris credit to Shimura. For a topological space $\mathcal{X}$ with a base point, addition of the base point defines a map $\Sym^n \mathcal{X} \to \Sym^{n+1} \mathcal{X}$, and the infinite symmetric product $\Sym^{\infty} \mathcal{X}$ is defined to be the direct limit $\Sym^{\infty} \mathcal{X} = \varinjlim_n \Sym^n \mathcal{X}$. Taking the union of two finite sets of points of $\mathcal{X}$ determines an addition on $\Sym^{\infty} \mathcal{X}$ whose identity is the base point.

\tpoint{Proposition}\label{pi_0SymXG_freeVS}{\em Let $\mathcal{X}$ be a path connected, Hausdorff, topological space with an action of $G=\Z/2$ equipped with a $G$-fixed base point. Assume that the path components of $\mathcal{X}^G$ are closed in $\mathcal{X}^G$. Then the monoid structure on $\pi_0((\Sym^{\infty} \mathcal{X})^G)$ determined by the monoid structure on $\Sym^{\infty} \mathcal{X}$ is a $\Z/2$ vector space structure, and $\pi_0 ((-)^G)$ applied to $\mathcal{X} \rightarrow \Sym^{\infty} \mathcal{X}$ is canonically isomorphic to the unit of the adjunction $(\mathcal{V}, \mathcal{I})$ on the pointed set $\pi_0(\mathcal{X}^G)$.}

Since $\Pic^0 X^+$ is an abelian group scheme, $\alpha$ determines a map $\Sym^n X \to \Pic^0 X^+$ from the $n$th symmetric product of $X$ to its generalized Jacobian. The following Proposition is well-known, but a proof is provided in the appendix for completeness.

\tpoint{Proposition}\label{SymnPicn}{\em For $n$ sufficiently large, $\Sym^n X \rightarrow \Pic^0 X^+ $ is an affine bundle (projectivized vector bundle) for $X$ non-proper (respectively proper).}

\epoint{Remark} It follows that when the real points of $\Sym^n X$ and $ \Pic^0 X^+ $ are given the analytic topology, $\Sym^n X (\R) \to \Pic^0 X^+ (\R)$ is an affine bundle or projectivized vector bundle for $n$ sufficiently large, giving an alternate way to see \cite[Prop 3.2 (2) $n(W^d)= n(S^d X)$ for large $d$]{GrossHarris} and \cite[2.7.4 $\mathcal{O} (D') \simeq \mathcal{O} (D)$]{Broglia}\hidden{strictly speaking, this is for $\deg D$ large, but translating $D$ by multiples of a real point shows the general case}.

The natural map from the $n$th symmetric power of $X(\C)$ to $\Sym^n X (\C)$\hidden{induced from $X \times X \times... \times X \to \Sym^n X$} is a homeomorphism\hidden{the map is a bijection. It extends to the analogous map for the smooth compactification of $X$, which is a continuous bijection between compact Hausdorff spaces, whence a homeomorphism}, so both may be denoted $\Sym^n X (\C)$. Let $\Sym^{\infty} X(\C)= \varinjlim_n  \Sym^n X (\C)$ be the infinite symmetric product of the topological space $X(\C)$.  The maps $\Sym^n X \rightarrow \Pic^0 X^+ $ are compatible with the map $\Sym^n X \to \Sym^{n+1} X$ given by addition of the base point, defining a $G_{\R}$-equivariant map $\Sym^{\infty} X(\C) \to  \Pic^0 X^+(\C).$ 

\tpoint{Proposition}\label{pi0(Syminfty_to_Pic)GR_iso}{\em $\pi_0 ((\Sym^{\infty} X(\C))^{G_{\R}}) \to \pi_0 ( \Pic^0 X^+(\R))$ is a bijection.}

\begin{proof}
The natural map from $ \varinjlim_n  \Sym^n X (\R)$ to the fixed points of $\Sym^{\infty} X(\C)$ is a homeomorphism\hidden{pf: there is a canonical continuous map $\varinjlim_n  \Sym^n X (\R) \to (\Sym^{\infty} X(\C))^{G_\R}$, which is bijective by inspection. A closed subset of $\varinjlim_n  \Sym^n X (\R)$ is a set $A$ whose inverse image $A_n$ in any $\Sym^n X (\R)$ is closed. Since $A_n$ is closed in $\Sym^n X(\C)$, we have that $A$ is also closed in $\Sym^{\infty} X(\C)$.}. By Proposition \ref{SymnPicn}, $\pi_0(\Sym^n X (\R)) \to \pi_0(\Pic^0 X^+ (\R))$ is a bijection for sufficiently large $n$, whence $\pi_0 (\varinjlim_n  \Sym^n X (\R)) = \varinjlim_n  \pi_0 (\Sym^n X (\R)) \to \pi_0(\Pic^0 X^+ (\R))$\hidden{homology, whence $\pi_0$, commutes with direct limits} is a bijection. 
\end{proof}

\epoint{Remark}\label{relation_condition_normalization_components_not_pointless_to_condition_X(C)_is_K(pi,1)} The hypothesis of Proposition \ref{kappaab_free_vs} and the resulting hypothesis of Theorem \ref{2nilsectionconjecture} is different from that of the real section conjecture, which is that $X(\C)$ be a $\rK(\pi,1)$. The normalization $\tilde{X} \to X$ induces a continuous map $\pi: \tilde{X}(\C) \to X(\C)$, which factors through the quotient $\tilde{X}(\C) \to Y$ where $Y$ is obtained by identifying points of $\tilde{X}(\C)$ with equal images under $\pi$. The homeomorphism $Y \to X(\C)$ shows that $X(\C)$ is homotopy equivalent to the wedge of the connected components of $\tilde{X}(\C)$ and a certain number of circles $S^1$. Since a wedge of $K(\pi,1)$'s is a $K(\pi,1)$, one sees that $X(\C)$ is a $K(\pi,1)$ if and only if none of the connected components of $\tilde{X}(\C)$ are $\proj^1_{\C}$.

\epoint{Remark}\label{pi0(SymGtoSymhG)bij}  Proposition \ref{kappaab_free_vs} and Proposition \ref{pi_0SymXG_freeVS} show that $$\pi_0 ((\Sym^{\infty}X(\C))^{G_{\R}}) \to  \pi_0((\Sym^{\infty}X(\C))^{h G_{\R}}) \to \mathscr{S}_{\pi_1((\Sym^{\infty}X(\C))/G_{\R})} $$ is a bijection. The second map is injective by the Dold-Thom theorem and the spectral sequence $\rH^i(G_{\R}, \pi_j (\Sym^{\infty} X(\C))) \Rightarrow \pi_{j-i} (\Sym^{\infty} X(\C)) ^{h G}$ \cite[IX \S 4]{Bousfield_Kan}. Thus the map from the fixed points to the homotopy fixed points of $\Sym^{\infty}X(\C)$ is a bijection on $\pi_0$. Note that $\Sym^{\infty} X(\C)$ is not a finite complex, so Sullivan's conjecture does not apply.

\section{$2$-nilpotent obstruction}\label{2nilobstruction_section}

Recall that for an extension of profinite groups \begin{equation}\label{piGextn} 1 \to \pi \to \tilde{\pi} \to G \to 1,\end{equation} the conjugacy class of a section $s: G \to \tilde{\pi}$ refers to the set of sections of the form $$g \mapsto \gamma s(g) \gamma^{-1}$$ where $\gamma$ is in $\pi$. Pushing out \eqref{piGextn} by $\pi \to \pi/[\pi]_n$ gives the extension \begin{equation}\label{pi/[pi]_nGextn} 1 \to \pi/[\pi]_n \to \tilde{\pi}/[\pi]_n \to G \to 1.\end{equation} When \eqref{piGextn} is equipped with a splitting, the extensions \eqref{pi/[pi]_nGextn} inherit splittings, which induce bijections between $\rH^1(G, \pi/[\pi]_n)$ and the conjugacy classes of sections of \eqref{pi/[pi]_nGextn}. Given a section $s: G \to \tilde{\pi}/[\pi]_2$ of \eqref{pi/[pi]_nGextn} for $n=2$, there exists a section $\tilde{s}: G \to \tilde{\pi}/[\pi]_3$ of \eqref{pi/[pi]_nGextn} for $n=3$ such that the composition $G \stackrel{\tilde{s}}{\to} \tilde{\pi}/[\pi]_3 \to \tilde{\pi}/[\pi]_2$ is in the conjugacy class of $s$ if and only if the class of $s$ vanishes under $$\delta_2: \rH^1(G, \pi/[\pi]_2) \to \rH^2(G, [\pi]_2/[\pi]_3)$$ where $\delta_2$ is the boundary map in continuous group cohomology from the extension $$1 \to [\pi]_2/[\pi]_3 \to \pi/[\pi]_3 \to \pi/[\pi]_2 \to 1.$$ 

We show that $\Ker \delta_2 =   \operatorname{Image} \kappa^{\ab}$ (Theorem \ref{kerdelta2=Imkappaab}) for $G= G_{\R}$ and $\tilde{\pi}$ equal to the \'etale or orbifold fundamental group of $X$ or $X(\C)/ G_{\R}$, respectively, with $X$ a curve over $\R$ satisfying the hypotheses of Theorem \ref{2nilsectionconjecture}. This proves Theorem \ref{2nilsectionconjecture}, since $\kappa^{\ab}$ is injective by Proposition \ref{kappaab_free_vs}.

$\delta_2$ is quadratic with associated bilinear form given by the following cup product \cite[p 242]{Zarkhin}: let $$[-,-]: \pi/[\pi]_2 \otimes \pi/[\pi]_2 \to [\pi]_2/[\pi]_3$$ be defined $[\gamma, \delta] = \tilde{\gamma} \tilde{\delta} \tilde{\gamma}^{-1} \tilde{\delta}^{-1}$ where $\tilde{\gamma}$, $\tilde{\delta}$ in $\pi/[\pi]_3$ map to $\gamma$, $\delta$ in $\pi/[\pi]_2$, respectively. (Since the choices of $\tilde{\gamma}$ differ by an element of the center, this map is well-defined on $\pi/[\pi]_2 \times \pi/[\pi]_2$. Bilinearity follows from \cite[Thm 5.1 p. 290]{Magnus_Karrass_Solitar}.) The cup product $$\rH^1(G, \pi/[\pi]_2) \otimes \rH^1(G, \pi/[\pi]_2) \to \rH^2(G, \pi/[\pi]_2 \otimes \pi/[\pi]_2)$$ can be pushed forward by $[-,-]$ to give a map \begin{equation}\label{[-,-]cup}\rH^1(G, \pi/[\pi]_2) \otimes \rH^1(G, \pi/[\pi]_2) \to \rH^2(G, [\pi]_2/[\pi]_3).\end{equation}

\tpoint{Proposition {\em(Zarkhin)}}\label{delta2quadratic}{\em For all $x$,$y$ in $\rH^1(G, \pi/[\pi]_2)$, $$\delta_2(x+y)= \delta_2(x) + \delta_2(y) + [-,-]_* x \cup y.$$ }

\epoint{Remark} For $\tilde{\pi}$ the \'etale fundamental group of a smooth, based, algebraic curve over a field and $G$ the absolute Galois group, Jordan Ellenberg introduced and studied $\delta_2$ as an obstruction to rational points of the Jacobian lying on the curve \cite{Ellenberg_2_nil_quot_pi}.

We will show injectivity of the bilinear form $(x,y) \mapsto [-,-]_* x \cup y$, for which we need the following notation and lemmas.

For an abelian group $\mathcal{L}$, let $\mathcal{L} \wedge \mathcal{L}$ denote the quotient of $\mathcal{L} \otimes \mathcal{L}$ by the relation $\ell \otimes \ell = 0$ for all $\ell$ in $\mathcal{L}$. When $G=\Z/2$, the pairing $$\cup: \rH^1(G, \mathcal{L}) \otimes \rH^1(G, \mathcal{L})  \to \rH^2(G, \mathcal{L} \wedge \mathcal{L})$$ induced by the cup product satisfies $$\ell \cup \ell = 0$$ by a straight forward computation, giving a natural map $$ \rH^1(G, \mathcal{L}) \wedge \rH^1(G, \mathcal{L})  \to \rH^2(G, \mathcal{L} \wedge \mathcal{L}).$$

\tpoint{Lemma}\label{cupwedgeinjectiveZ/2}{\em Let $\mathcal{L}$ be an abelian group or profinite abelian group with no $2$-torsion and an action of $G=\Z/2$. Then the natural map $$\rH^1(G, \mathcal{L}) \wedge \rH^1(G, \mathcal{L}) \rightarrow \rH^2(G, \mathcal{L} \wedge \mathcal{L}) $$ induced by the cup product is injective.}

\begin{proof}
Let $\tau$ denote the generator of $G$. Since $G$ has order $2$, $\rH^1(G, \mathcal{L})$ is a $\Z/2$-vector space. The short exact sequence $$\xymatrix{0 \ar[r] & \mathcal{L} \ar[r]^{2} & \mathcal{L} \ar[r] & \mathcal{L} / 2 \mathcal{L} \ar[r] & 0} $$ shows that $$\rH^1(G, \mathcal{L}) \hookrightarrow \rH^1(G, \mathcal{L}/ 2 \mathcal{L})$$ is injective. Similarly, $$\rH^2(G, \mathcal{L} \wedge \mathcal{L}) \hookrightarrow \rH^2(G, \mathcal{L}/2 \wedge \mathcal{L}/2)$$ is injective. Thus, it suffices to show that $$\rH^1(G, W) \wedge \rH^1(G, W) \rightarrow \rH^2(G, W \wedge W) $$ is injective for any $\Z/2$-vector space $W$ with a $G$ action. This map is described by $$ [w_1] \wedge [w_2] \mapsto [w_1 \wedge \tau w_2] $$ where $w_i$ is in the kernel of $\tau + 1: W \rightarrow W$, and $[w_i]$ denotes the corresponding cohomology class via the cyclic $G$ resolution of $\Z$ \cite[V \S 1 pg. 108]{Brown_coh_groups}.

Note that $(1+ \tau) (w_1 \wedge w_2) = w_1  \wedge w_2 + \tau w_1 \wedge \tau w_2 = (w_1 + \tau w_1)  \wedge w_2 + \tau w_1 \wedge (w_2 +\tau w_2) $. Therefore, we have a map $H^2(G, W \wedge W) \rightarrow (W/I) \wedge (W/I)$ where $I$ denotes the image of $\tau +1:W \rightarrow W$.

Let $K$ denote the kernel of $\tau + 1 : W \rightarrow W$. Since the automorphisms $\tau +1$ and $\tau -1$ of $W$ are equal, $\rH^1(G,W)\cong K/I$. The inclusion $K \hookrightarrow W$ induces an injection $K/I \wedge K/I \hookrightarrow W/I \wedge W/I$ because $W$ is a vector space. 

The the desired injectivity follows from the commutative diagram $$\xymatrix{ \rH^1(G, W) \wedge \rH^1(G, W) \ar[rr]^{\cup} \ar[d]^{\cong} && \rH^2(G, W \wedge W) \ar[d] \\ K/I \wedge K/I \ar[rr] && (W/I) \wedge (W/I)}.$$
\end{proof}

Note that $[-,-]$ factors through $ \pi/[\pi]_2 \wedge \pi/[\pi]_2$ giving a map, also denoted $[-,-]$, $$[-,-]: \pi/[\pi]_2 \wedge \pi/[\pi]_2 \to [\pi]_2/[\pi]_3.$$

\tpoint{Lemma}\label{[-,-]*injective}{\em Let $X$ be a geometrically connected, based curve over $\R$ and let $\pi$ denote the \'etale or topological fundamental group of $X_{\C}$, $X(\C)$ respectively. Then $[-,-]_*: \rH^2(G_{\R}, \pi/[\pi]_2 \wedge \pi/[\pi]_2) \to \rH^2(G_{\R}, [\pi]_2/[\pi]_3)$ is injective.}

\begin{proof}

By \cite[Lemma 1.3]{Dwyer}, for any group $\pi$, there is a right exact sequence \begin{equation}\label{Dwyereq}\rH_2(\pi) \to \rH_2(\pi/[\pi]_2) \to [\pi]_2/[\pi]_3 \to 1, \end{equation} where $\rH_2$ denotes group homology with $\Z$-coefficients. When $\pi/[\pi]_2$ is a free $\Z$-module, $\rH_2(\pi/[\pi]_2)$ is canonically identified with $\pi/[\pi]_2 \wedge \pi/[\pi]_2$ and the surjection of \eqref{Dwyereq} is $[-,-]$.

Let $\pi = \pi_1^{\topp}X(\C)$, so $\pi/[\pi]_2$ is a free $\Z$-module. Since complex conjugation is orientation reversing, the canonical surjection $\rH_2(X(\C)) \to \rH_2(\pi)$ \cite[Thm 5.2 p. 41]{Brown_coh_groups} shows that $G_{\R}$ acts on $\rH_2(\pi)$ by multiplication by $-1$. Thus $G_{\R}$ acts on the kernel $K$ of $[-,-]: \pi/[\pi]_2 \wedge \pi/[\pi]_2 \to [\pi]_2/[\pi]_3$ by multiplication by $-1$. Since $\pi/[\pi]_2 \wedge \pi/[\pi]_2$ is free, so is $K$. Thus $\rH_2(G_{\R}, K)=0$, giving the desired injection in this case.

For $\mathcal{L}$ a free $\Z$-module with $\Z/2$-action, $\rH^i(\Z/2, \mathcal{L}) \to \rH^i(\Z/2, \mathcal{L}^{\wedge})$ induced by profinite completion $ \mathcal{L} \to \mathcal{L}^{\wedge}\cong \mathcal{L} \otimes_{Z} \Zhat$ is an isomorphism because $\Zhat$ is torsion-free whence flat over $\Z$. Since $\pi/[\pi]_2$ and $[\pi]_2/[\pi]_3$ are free $\Z$-modules, the \'etale case follows. 
\end{proof}

Let $X$ be as in Theorem \ref{2nilsectionconjecture} and define $\delta_2$ with the extension \eqref{htpyes} or its topological analogue.

\tpoint{Theorem}\label{kerdelta2=Imkappaab}{\em $\Ker \delta_2 = \operatorname{Image} \kappa^{\ab}$.}

\begin{proof}

For $p$ in $X(\R)$, $\kappa(p)$ determines a section of \eqref{pi/[pi]_nGextn} for $n=3$, so $\delta_2(\kappa^{\ab}(p)) = 0$.

By Proposition \ref{kappaab_free_vs}, an arbitrary element of $\rH^1(G_{\R}, \pi^{ab})$ is of the form $x_1 + x_2 + \ldots + x_n$ where the $x_i$ are in the image of $\kappa^{\ab}$ and $\{ x_1, x_2,\ldots, x_n \}$ is linearly independent. By Proposition \ref{delta2quadratic}, $$\delta_2(x_1 + x_2 + \ldots + x_n) = \sum_{i=1}^n \delta_2(x_i) + \sum_{1\leq i < j \leq n} [-,-]_* x_i \cup x_j.$$ Furthermore, $\delta_2 (x_i) = 0$ since $x_i \in\operatorname{Image} \kappa^{\ab}$. If $n>1$, then $ \sum_{1\leq i < j \leq n} x_i \wedge x_j$ is a non-zero element of $\rH^1(G_{\R}, \pi^{\ab}) \wedge \rH^1(G_{\R}, \pi^{\ab})$, whence $\sum_{1 \leq i < j \leq n} [-,-]_* x_i \cup x_j$ is a non-zero element of $\rH^2(G_{\R}, [\pi]_2/[\pi]_3)$ by Lemmas \ref{cupwedgeinjectiveZ/2} and \ref{[-,-]*injective}. Thus $n=1$ for any element of $\Ker \delta_2 $, showing $\Ker \delta_2 \subset  \operatorname{Image} \kappa^{\ab}$. \end{proof}

\section{Birational $2$-nilpotent real section conjecture}\label{birational_section}

Let $X$ be a smooth, proper, geometrically connected curve over $\R$ such that $X(\R) \neq \emptyset$. Let $X(\R)^{\pm}$ denote the set of real points of $X$ equipped with a real tangent direction. 

Choose a local parameter $z \in \R(X)$ at a point of $X(\R)$, and note that $z$ embeds the function field $\R(X)$ into the field of Puiseux series \[
\C((z^\Q)) : = \cup_{n\in \Z_{>0}} \C ((z^{1/n})),
\] which is algebraically closed. Taking the algebraic closure $\Omega_z$ of $\R(X)$ in $\C((z^\Q))$ and considering the intermediate extension $\R(X) \subset \C(X)$ gives \begin{equation}\label{Gal_hes} 1 \to \Gal(\Omega_z/\C(X)) \to \Gal(\Omega_z/\R(X)) \to G_{\R} \to 1 \end{equation} which should be viewed as an analogue of the homotopy exact sequence \eqref{htpyes}. The coefficientwise action of $G_{\R}$ on  $\C((z^\Q))$ defines a section of \eqref{Gal_hes}. 

Given a second local parameter $w$ and the associated embeddings $$\R(X) \subset \overline{\R(X)} =: \Omega_w \subset \C ((w^\Q))$$ choose an isomorphism $\Omega_z \cong \Omega_w$ which is the identity on the inclusion of $\C(X)$ in both fields, determining an isomorphism $\Gal(\Omega_z/\R(X)) \cong \Gal(\Omega_w/\R(X))$. The map $G_{\R} \to \Gal(\Omega_w/\R(X))$ defined by the coefficientwise action of $G_{\R}$ on  $\C((w^\Q))$ then gives a second section of \eqref{Gal_hes}. It is not difficult to check that the conjugacy class of this section depends only on the element of $X(\R)^{\pm}$ determined by the tangent vector $$\Spec \R[w]/\langle w^2 \rangle \to X $$ associated to $w$.

It follows that $b \in X(\R)^{\pm}$ determines a map $\kappa: X(\R)^{\pm} \to \mathscr{S}_{\Gal(X/\R)},$ where $\mathscr{S}_{\Gal(X/\R)}$ denotes the set of conjugacy classes of sections of $$1 \to G_{\C(X)} \to G_{\R(X)} \to G_{\R} \to 1.$$ 

We introduce some notation. Let $\mathscr{S}_{\Gal(X/\R)}^n$ denote the set of conjugacy classes of sections of the push-out sequence \begin{equation}\label{n_Gal_hes} 1 \to G_{\C(X)}/[G_{\C(X)}]_n \to G_{\R(X)} /[G_{\C(X)}]_n \to G_{\R} \to 1,\end{equation} and for $m > k$, let $\mathscr{S}_{\Gal(X/\R)}^{k \leftarrow m}$ denote the set of conjugacy classes of sections of \eqref{n_Gal_hes} for $n=k$ which have a lift to a section of \eqref{n_Gal_hes} for $n=m$. For $V \subseteq X$ Zariski open, let $\mathscr{S}_{\pi_1(V/\R)}^{n}$ denote the set of conjugacy classes of sections of \begin{equation}\label{n_pi_hes} 1 \to \pi_1(V_{\C})/[\pi_1(V_{\C})]_n \to \pi_1(V)/[\pi_1(V_{\C})]_n \to G_{\R} \to 1,\end{equation} where $\pi_1$ denotes the \'etale fundamental group, and let $\mathscr{S}_{\pi_1(V/\R)}^{k \leftarrow m}$ denote the conjugacy classes of sections of \eqref{n_pi_hes} for $n=k$ which have a lift to a section of \eqref{n_pi_hes} for $n=m$.

\tpoint{Corollary}\label{corbir2nilSC}{\em $\kappa$ gives a natural bijection from $X(\R)^{\pm}$ to $\mathscr{S}_{\Gal(X/\R)}^{2 \leftarrow 3}$.}

\begin{proof}
For $U$ a Zariski open of $X$, applying $\pi_1$ to the inclusion of the generic point gives a map $G_{\R(X)} \to \pi_1(U,b)$ \cite[V Proposition 8.1]{sga1}. By functoriality of $\pi_1$,\hidden{
For any $V \subset U \subset X$, we have $$\xymatrix{ 1 \ar[r] & G_{\C(X)}/[G_{\C(X)}]_n \ar[d] \ar[r] & G_{\R(X)}/[G_{\C(X)}]_n \ar[r] \ar[d] & G_{\R} \ar[d] \ar[r] & 1\\ 
1 \ar[r] & \pi_1(V_{\C})/[\pi_1(V_{\C})]_n \ar[r] \ar[d] & \pi_1(V)/[\pi_1(V_{\C})]_n \ar[r] \ar[d] & G_{\R} \ar[r] \ar[d] & 1 \\
1 \ar[r] & \pi_1(U_{\C})/[\pi_1(U_{\C})]_n \ar[r] & \pi_1(U)/[\pi_1(U_{\C})]_n \ar[r] & G_{\R} \ar[r] & 1 } $$ where all fundamental groups are based at $b$. 
} for any $V \subset U \subset X$, we have compatible maps $$\mathscr{S}_{\Gal(X/\R)}^{n} \to \mathscr{S}_{\pi_1(V/\R)}^{n}  \to \mathscr{S}_{\pi_1(U/\R)}^{n}$$  $$\mathscr{S}_{\Gal(X/\R)}^{n \leftarrow m} \to \mathscr{S}_{\pi_1(V/\R)}^{n \leftarrow m}  \to \mathscr{S}_{\pi_1(U/\R)}^{n \leftarrow m}$$ for any $n$ and $m > n$. These maps are compatible with $\kappa$ in that the diagram $$\xymatrix{ X(\R)^{\pm}  \ar[d] \ar[r] & \ar[d] \mathscr{S}_{\Gal(X/\R)}^{2 \leftarrow 3} \\
\varprojlim_{U} \pi_0(U(\R)) \ar[r] & \varprojlim_U \mathscr{S}_{\pi_1(U/\R)}^{2 \leftarrow 3}} $$ commutes, where $X(\R)^{\pm} \to  \pi_0(U(\R))$ sends an element of $X(\R)^{\pm}$ to the connected component the tangent direction points to, and where $U$ runs over the Zariski opens of $X$.

By Theorem \ref{2nilsectionconjecture}, $\kappa: \pi_0(U(\R)) \to \mathscr{S}_{\pi_1(U/\R)}^{2 \leftarrow 3}$ is a bijection, showing that the bottom horizontal arrow is a bijection. The left vertical arrow is a bijection by inspection. It thus suffices to show that $\mathscr{S}_{\Gal(X/\R)}^{2} \to \varprojlim_U \mathscr{S}_{\pi_1(U/\R)}^{2}$ is injective, which is equivalent to showing that $\rH^1(G_{\R}, G_{\C(X)}^{\ab}) \to \varprojlim_U \rH^1(G_{\R}, \pi_1(U_{\C})^{\ab})$ is injective \cite[IV 2.3]{Brown_coh_groups}.

By \cite[V Proposition 8.2]{sga1}\hidden{the referenced Proposition gives $\pi_1$ of the normal schemes $U$ as the Galois group of the maximal unramified extension. To complete the proof, one needs to know that for any finite $L$ extension of $\C(X)$ there is a Zariski open $U \subset X$ over which $L$ is unramified. Let $X_L \to X$ denote the normalization of $X$ in $L$. We must show that the sheaf of relative differentials $\Omega_{X_L/X}$ is supported on a proper closed subscheme of $X_L$. Because the stalk of $\Omega_{X_L/X}$ at the generic point is $\Omega_{L/\C(X)} = 0$, we have that the stalk of $\Omega_{X_L/X}$ at the generic point vanishes. Since $\Omega_{X_L/X}$ is a coherent sheaf, we $\Omega_{X_L/X}$ is supported on a proper closed subscheme of $X_L$} the natural map $G_{\C(X)} \to \varprojlim_{U} \pi_1(U_{\C}, b)$ is an isomorphism. It follows that $G_{\C(X)}^{\ab} \to \varprojlim_{U} \pi_1(U_{\C}, b)^{\ab}$ is an isomorphism.\hidden{

Let $G$ be a profinite group. Then by RZ 2.1.3, $G = \varprojlim_{i \in I} G_i$ with $G_i = G/U_i$, $i \in I$, where $I$ is a directed set and $U_i$ are open normal subgroups of $G$. We claim that $G^{\ab} \to \varprojlim_{i \in I} G_i^{ab}$ is an isomorphism, where $G^{\ab}= G/ \overline{[G,G]}$:

$G \to G_i$ is a surjection by Math 131 PS62(b). Thus we have surjections $G \to G_i^{\ab}$ for all $i$. Thus we have that $ G \to \varprojlim_{i \in I} G_i^{\ab}$ is a surjection because given $(g_i)_{i \in I} \in  \varprojlim_{i \in I} G_i^{\ab}$,  the subset $X_i$ of those elements of $G$ mapping to $g_i$ is a closed non-empty set, and the collection $\{ X_i \}$ has the finite intersection property. Thus $G^{\ab} \to \varprojlim_{i \in I} G_i^{\ab}$ is a surjection. Take $g \in G$ mapping to a non-zero element of $G^{\ab}$. Then there is an open normal subgroup $U$ of $G$ containing $\overline{[G,G]}$ and not containing $g$ by RZ 2.1.4(d). By RZ Prop 2.1.5(a), there exists $U_i$ such that $U_i$ is a subgroup of $U$. Thus we have $G \to G/U_i \to G/U.$ Since $U$ contains $\overline{[G,G]}$, $G/U$ is an abelian group, whence we have $G \to (G/U_i)^{\ab} \to G/U.$ Since the image of $g$ is non-zero in $G/U$, the image of $g$ is non-zero in $(G/U_i)^{\ab}$. Thus $G^{\ab} \to \varprojlim_{i \in I} G_i^{\ab}$ is injective.} Since $\varprojlim$ is exact as a functor on inverse systems of compact abelian groups\hidden{this is stated without proof pg 153 Prop 3.2 Helmut Koch ``Algebraic number theory."}, the natural map $$\rH^1(G_{\R}, \varprojlim_U \pi_1(U_{\C})^{\ab}) \to \varprojlim_U \rH^1(G_{\R}, \pi_1(U_{\C})^{\ab})$$ is an isomorphism, showing the corollary.\hidden{
\tpoint{Lemma}\label{H1lim=limH1}{\em Let $A$ be a functor from the opposite category of a directed set $I$ to profinite abelian groups with an action of $G=\Z/2$, and where morphisms are required to be surjective. Then the natural map $\rH^1(G, \varprojlim_I A) \to \varprojlim_I \rH^1(G, A)$ is an isomorphism.}

\begin{proof}
Let $U$ denote an object of $I$.

$$0 \to \Image(1 - \tau: A(U)) \to \Ker (1 + \tau : A(U)) \to \rH^1(G_{\R}, A(U)) \to 0.$$

Taking $\varprojlim$ gives the exact sequence 

$$0 \to \varprojlim_U \Image(1 - \tau : A(U)) \to \varprojlim_U \Ker (1 + \tau : A(U) )\to \varprojlim_U \rH^1(G_{\R}, A(U)) \to 0.$$

$ \varprojlim_U \Ker (1 + \tau : A(U) )= \Ker (1 + \tau : \varprojlim_U  A(U))$.

It is thus sufficient to show that the natural map $ \Image(1 - \tau : \varprojlim_U A(U)) \to \varprojlim_U \Image(1 - \tau : A(U) )$ is a surjection.

$0 \to \Ker (1 - \tau: A(U)) \to  A(U) \to \Image(1 - \tau : A(U) )\to 0.$

We need the sequence obtained by applying $\varprojlim$ to be right exact: take $(a_U)_{U}$ in $\varprojlim \Image(1 - \tau) A(U)$. Let $X(U)$ be the subset of elements of $\varprojlim A(U)$ whose image under $1- \tau$ has $U$-coordinate equal to $a(U)$. $X(U)$ is a non-empty closed subset of $\varprojlim A(U)$. The collection $X(U)$ has the finite intersection property because for all $U_1$, $U_2$,\ldots $U_n$, $X(U_1 \cap \ldots U_n) \subset \cap_{i=1}^n X(U_i)$. Since $\varprojlim A(U)$ is compact, it follows that $\cap_U X(U)$ is non-empty, showing that $\varprojlim A(U) \to \varprojlim \Image(1 - \tau) A(U)$ is surjective as desired.  

\end{proof}}
\end{proof}

\section{$2$-nilpotent topological approximation}\label{top2nilapprox}

Let $X$ be a geometrically connected curve over $\R$ equipped with a base point in $X(\R)$, and let $\pi$ denote the topological fundamental group of $X(\C)$.

\epoint{Abelian approximation} The action of $\pi /[\pi]_2$ on itself by left translation gives an injective homomorphism $$\pi /[\pi]_2 \to \operatorname{Affine}( \pi /[\pi]_2),$$ where $\operatorname{Affine}(\pi /[\pi]_2)$ denotes the group of invertible affine transformations of the free abelian group $\pi/[\pi]_2$. (The image of $\pi /[\pi]_2$ is contained in the subgroup of translations, but\hidden{ this will change when $\pi /[\pi]_2$ is replaced by $\pi /[\pi]_3$, so} the notation is set up for the larger group.) Tensoring with $\R$ gives an action of $\pi /[\pi]_2$ on $\pi /[\pi]_2 \otimes_{\Z} \R$, and taking the quotient gives a model for $\rK(\pi /[\pi]_2, 1)$ denoted $$\Alb_1 = (\pi/[\pi]_2) \backslash (\pi /[\pi]_2 \otimes_{\Z} \R).$$ 

The Galois group $G_{\R}$ acts on $\pi/[\pi]_2$, giving a linear action on $\pi /[\pi]_2 \otimes_{\Z} \R$. For $g \in G_{\R}$ and $\gamma \in \pi/[\pi]_2$, we have an element $g \gamma $ of $\pi/[\pi]_2$, and the equality $g(\gamma v) = (g \gamma) (gv)$ for all $v \in \pi /[\pi]_2 \otimes_{\Z} \R$. Thus $\Alb_1$ inherits a $G_{\R}$-action.

\epoint{$2$-Nilpotent approximation} Choose a set of generators $x_1, \ldots, x_n$ for $\pi$ which are a basis for $\pi/[\pi]_2$.\hidden{This is possible because $X(\C)$ is homotopy equivalent to a wedge of $S^1$'s and orientable surfaces} Let $s$ be the set-theoretic section of the quotient map $q: \pi/[\pi]_3 \to \pi/[\pi]_2$ given by taking $x_1^{a_1} x_2^{a_2} \cdots x_n^{a_n}$ in $\pi/[\pi]_2$ to  $x_1^{a_1} x_2^{a_2} \cdots x_n^{a_n}$ in $\pi/[\pi]_3$. The section $s$ determines a bijection between $\pi/[\pi]_2 \oplus [\pi]_2/[\pi]_3$ and $\pi/[\pi]_3$ by sending $v \oplus z \in \pi/[\pi]_2 \oplus [\pi]_2/[\pi]_3$ to $s(v) z$. Via this bijection, the group law on $\pi/[\pi]_3$ gives a composition law on $\pi/[\pi]_2 \oplus [\pi]_2/[\pi]_3$, denoted $\circ$. Let $+$ denote addition on the free abelian group $\pi/[\pi]_2 \oplus [\pi]_2/[\pi]_3$. There is a bilinear pairing $$\langle - , - \rangle: \pi/[\pi]_2 \oplus \pi/[\pi]_2 \to  [\pi]_2/[\pi]_3$$ such that $$v \circ w = v + w + \langle q v, q w \rangle,$$ for $v$ and $w$ in $\pi/[\pi]_2 \oplus  [\pi]_2/[\pi]_3$. Here and later, we slightly abuse the notation $q$ by letting $q$ also denote the projection $\pi/[\pi]_2 \oplus  [\pi]_2/[\pi]_3 \to \pi/[\pi]_2$. Thus the action of $\pi /[\pi]_3$ on itself by left translation gives an injective homomorphism $$\pi /[\pi]_3 \to \operatorname{Affine}(\pi/[\pi]_2 \oplus [\pi]_2/[\pi]_3) \subset \operatorname{Affine}( (\pi/[\pi]_2 \oplus [\pi]_2/[\pi]_3) \otimes_{\Z} \R).$$  For example, $x_2$ is sent to the affine transformation \begin{equation*}(a_1,a_2,\ldots,a_n) \times z  \mapsto (a_1,a_2 + 1,\ldots,a_{n}) \times (z - a_1 [x_1,x_2]).\end{equation*}

Taking the quotient of $(\pi/[\pi]_2 \oplus [\pi]_2/[\pi]_3) \otimes_{\Z} \R$ by $\pi/[\pi]_3$ gives a model for $\rK(\pi /[\pi]_3, 1)$, denoted $$\Alb_2 = (\pi/[\pi]_3) \backslash ((\pi/[\pi]_2 \oplus [\pi]_2/[\pi]_3) \otimes_{\Z} \R).$$

$G_{\R}$ acts on $\pi/[\pi]_3$. Give $\pi/[\pi]_2 \oplus [\pi]_2/[\pi]_3$ the $G_{\R}$-action inherited from its bijection with $\pi/[\pi]_3$. Since for all $z \in 0 \oplus [\pi]_2/[\pi]_3$ we have $z \circ w = w \circ z = w + z$, it follows that $$g (v \circ w) = g( (v + w) \circ \langle q v, q w \rangle ) = g(v + w) \circ g \langle q v, q w \rangle = g(v + w) + g \langle q v, q w \rangle,$$ $$(gv) \circ (gw) = (gv) + (gw) + \langle q g v, q gw \rangle.$$ Since $g (v \circ w) = (g v) \circ (gw)$ for $g \in G_{\R}$, we have that $$g (w + v) - g(w) - g(v) = \langle g q  v, g q w \rangle - g \langle q v, q w \rangle$$ is bilinear. Let $B_g(v,w) = \langle g q  v, g q w \rangle - g \langle q v, q w \rangle $ denote this bilinear form. It follows that \begin{equation}\label{GRAlb3action}g (v) = \frac{1}{2} B_g(v,v) + L_g(v)\end{equation} where $L_g$ is a linear endomorphism of $(\pi/[\pi]_2 \oplus [\pi]_2/[\pi]_3)\otimes_{\Z} \Z[\frac{1}{2}]$.\hidden{ the action of $G_{\R}$ on $\pi/[\pi]_2 \oplus [\pi]_2/[\pi]_3$ can be expressed as an ordered tuple of polynomials of degree $\leq 2$.} In particular, the $G_{\R}$ action on $\pi/[\pi]_2 \oplus [\pi]_2/[\pi]_3$ extends to $(\pi/[\pi]_2 \oplus [\pi]_2/[\pi]_3) \otimes_{\Z} \R.$ As above, the equality $g(\gamma v) = (g \gamma) (gv)$ for all $v \in (\pi/[\pi]_2 \oplus [\pi]_2/[\pi]_3) \otimes_{\Z} \R$, $g \in G_{\R}$ and $\gamma \in \pi/[\pi]_3$, implies that $G_{\R}$ acts on $\Alb_2$. 

The map $\Alb_2 \to \Alb_1$ induced by the projection $q$ is a $G_{\R}$-equivariant fiber bundle with fiber $\rK([\pi]_2/[\pi]_3,1).$ 

\epoint{Mapping $X(\C)$ to its approximations} Equip $\Alb_2$ with the base point induced by the origin of $(\pi/[\pi]_2 \oplus [\pi]_2/[\pi]_3) \otimes_{\Z} \R$, and say that a map between spaces with base points is pointed if the image of the base point of the domain is the base point of the codomain. Choose a pointed map $f: X(\C) \to \Alb_2$ such that the induced map $f_*$ on fundamental groups is the quotient. Let $\widetilde{X(\C)} \to X(\C)$ and $\widetilde{\Alb_2} \to \Alb_2$ denote the universal pointed covering maps of $X(\C)$ and $\Alb_2$ respectively\hidden{, which are unique up to unique isomorphism}. Note that $\widetilde{\Alb_2} \to \Alb_2$ is $(\pi/[\pi]_2 \oplus [\pi]_2/[\pi]_3) \otimes_{\Z} \R \to \Alb_2$, up to unique pointed isomorphism, where the origin is the base point of $(\pi/[\pi]_2 \oplus [\pi]_2/[\pi]_3) \otimes_{\Z} \R$. There is a unique lift of $f$ to a pointed map $\tilde{f}: \widetilde{X(\C)} \to \widetilde{\Alb_2}$. For all $\gamma \in \pi$, we can view $\gamma$ as an automorphism of $\widetilde{X(\C)}$ or $\widetilde{\Alb_2}$, and the choice of $f_*$ implies that $\gamma \tilde{f} = \tilde{f} \gamma$. Give $\widetilde{X(\C)}$ and $\widetilde{\Alb_2}$ the $G_{\R}$-actions lifting those of $X(\C)$ and $\Alb_2$ and such that the base points are fixed under $G_{\R}$. This $G_{\R}$-action on $(\pi/[\pi]_2 \oplus [\pi]_2/[\pi]_3) \otimes_{\Z} \R$ is consistent with the one above. Let $\tau$ denote complex conjugation and note that $\tau f \tau^{-1} : X(\C) \to \Alb_2$ is a pointed map such that the induced map on $\pi_1$ is the quotient. $\tau \tilde{f} \tau^{-1}$ is the unique lift of $\tau f \tau^{-1}$ to a pointed map between the pointed universal covering spaces, and $\gamma \tau \tilde{f} \tau^{-1} = \tau \tilde{f} \tau^{-1} \gamma$. 

Define $\tilde{g}: \widetilde{X(\C)} \to (\pi/[\pi]_2 \oplus [\pi]_2/[\pi]_3) \otimes_{\Z} \R = \widetilde{\Alb_2}$ by $$ \tilde{g} = \frac{1}{2} \tilde{f} + \frac{1}{2} \tau \tilde{f} \tau^{-1}.$$ By the above, $\gamma$ acts by affine transformations on $(\pi/[\pi]_2 \oplus [\pi]_2/[\pi]_3) \otimes_{\Z} \R$, which implies that $\gamma \tilde{g} = \frac{1}{2} \gamma \tilde{f} + \frac{1}{2} \gamma \tau \tilde{f} \tau^{-1} = \frac{1}{2} \tilde{f} \gamma+ \frac{1}{2}  \tau \tilde{f} \tau^{-1} \gamma = \tilde{g} \gamma$. Thus, $\tilde{g}$ induces a pointed map $g: X(\C) \to \Alb_2$ such that $g_*$ the quotient map on fundamental groups. Since $\tau$ acts linearly on the universal cover $\pi/[\pi]_2 \otimes_{\Z} \R$ of $\Alb_1$, we have that $q \tilde{g} : \widetilde{X(\C)} \to \pi/[\pi]_2 \otimes_{\Z} \R$ is $G_{R}$-equivariant. Thus for any $\tilde{x}$ in $\widetilde{X(\C)}$, the points $\tilde{g}$ and $\tau \tilde{g} \tau^{-1} $ are contained in the same fiber of $q$. The calculation \eqref{GRAlb3action} implies that $\tau$ determines an affine transformation between the fibers of $q: \widetilde{\Alb_2} \to \widetilde{\Alb_1}$ over $q \tilde{g}(\tilde{x})$ and $\tau q \tilde{g}(\tilde{x}) = q \tilde{g}(\tau \tilde{x})$. 

Define $\tilde{h}: \widetilde{X(\C)} \to \widetilde{\Alb_2}$ by $\tilde{h} = \frac{1}{2} \tilde{g} + \frac{1}{2} \tau \tilde{g} \tau^{-1}$. As above, $\tilde{h}$ induces a pointed map $h: X(\C) \to \Alb_2$ such that $h_*$ the quotient map on fundamental groups. Furthermore, $$\tau \tilde{h} \tau^{-1} = \tau(\frac{1}{2} \tilde{g} \tau^{-1} + \frac{1}{2} \tau \tilde{g}) = \frac{1}{2} \tau \tilde{g} \tau^{-1} + \frac{1}{2} \tilde{g} = \tilde{h}$$ where the second to last equality follows because $\tau$ is affine on the fiber over $q \tilde{g}(\tau^{-1} \tilde{x})$ for all $\tilde{x}$ in $\widetilde{X(\C)}$. Thus $\tilde{h}$ and $h$ and $G_{\R}$-equivariant. Use the notation $\alpha_2$ for $h$, so $\alpha_2 = h$. Let $\alpha: X(\C) \to \Alb_1$ denote the composition of $\alpha_2$ with the projection. The notation is chosen to recall that for the Abel-Jacobi map. 

\hidden{

Let $X$ and $Y$ be $G$-CW complexes with chosen $G$-fixed base points and suppose that $Y$ is a $\rK(\pi,1)$. A $G$-equivariant group homomorphism $\pi_1(X) \to \pi$ is not necessarily realizable as $\pi_1$ applied to a $G$-equivariant pointed map $X \to Y$. For example, let $Y = \R \proj^{\infty}$ with the trivial $G=\Z/2$ action, and let $X$ be $S^1$ viewed as the unit circle in $\C$ equipped with the action of complex conjugation. Then the non-trivial map $\pi_1(X) \cong \Z \to \pi_1 (\R \proj^{\infty}) \cong \Z/2$ is $G$-equivariant. However, any $G$-equivariant map from $S^1$ to a space with a trivial $G$-action factors through $S^1/ G = [0,1]$ and is therefore null-homotopic.

{\bf Question:} Given a finite $G$-CW complex $X$ with a chosen fixed point and a $\rK(\pi, 1)$ with a $G$-action and a chosen fixed point, what is the obstruction to the existence of a $G$-equivariant map $X \to \rK(\pi, 1)$ realizing a given $G$-equivariant group homomorphism $\pi_1(X) \to \pi$? How can such maps be classified up to $G$-equivariant homotopy? 

{\bf Partial answer:} The obstruction theory and classification should be the same as for CW complexes, but one would first apply this theory to the fixed points. Then to cells with some maximal subgroup as a stabilizer, etc. For example, for $\Z/2$, one would do obstruction theory first on the fixed points, and then on the free cells.

}

\epoint{Remark}\label{X_Alb_n_finite_GCW} $X(\C)$, $\Alb_1$, and $\Alb_2$ can be given the structure of finite $G_{\R}$-CW complexes in the sense of \cite{Bredon}; by the main theorem of \cite{Illman}, a smooth compact manifold with a $\Z/2$ action can be given such a structure, showing the existence of finite $G_{\R}$-CW complex structures for $\Alb_1$, $\Alb_2$, and the complex points of the smooth compactification of the normalization of X. Removing and or identifying finitely many points of a finite $G_{\R}$-CW complex yields a finite $G_{\R}$-CW complex\hidden{remove any simplices touching removed points}, so $X(\C)$ also has the structure of a finite $G_{\R}$-CW complex.

We obtain the commutative diagram of $G_{\R}$-equivariant maps between finite $G_{\R}$-CW complexes \begin{equation}\label{XAlb2Alb3CD}\xymatrix{ & \Alb_2 \ar[d] \\
X(\C) \ar[ur]^{\alpha_2} \ar[r]^{\alpha} & \Alb_1} \end{equation} and view $\Alb_1$ as an abelian approximation to $X$ and $\Alb_2$ as a $2$-nilpotent approximation to $X$. For $X$ smooth and proper, integration gives a natural map from $\pi/[\pi]_2$ to the $\C$-linear dual of the global holomorphic one-forms on $X$, denoted $\rH^0(X, \Omega)^*$. This map extends to a $G_{\R}$-equivariant $\R$-linear isomorphism $$\widetilde{\Alb_1}= \pi/[\pi]_2 \otimes_{\Z} \R \to \rH^0(X, \Omega)^*$$ which identifies $\Alb_1$ with the $G_{\R}$-topological space underlying the complex points of the Albanese variety of $X$. Thus $\Alb_1$ in \eqref{XAlb2Alb3CD} can be replaced with $\Pic^0 X (\C)$. The notation $\Alb$ comes from the view point that a $2$-nilpotent approximation to $X$ is analogous to a higher Albanese variety. See \cite{Hain_Higher_Albanese}, \cite{HainZucker}.

Theorem \ref{2nilsectionconjecture} can be rephrased as the statement that the connected components of real points of the curve are those of the Albanese which can be lifted to the $2$-nilpotent approximation:

\tpoint{Theorem}\label{Top_2nilSC}{\em Let $X$ be a geometrically connected, based curve over $\R$, such that each irreducible component of its normalization has $\R$-points. Let $\Alb_1$, $\Alb_2$, and $\alpha$ be as constructed above to obtain \eqref{XAlb2Alb3CD}. Then $\alpha$ induces a bijection from $\pi_0(X(\R))$ to the image of $ \pi_0(\Alb_2^{G_{\R}}) \to \pi_0(\Alb_1^{G_{\R}})$.}

\begin{proof}
By Remark \ref{X_Alb_n_finite_GCW}, $\Alb_1$ is a finite $G_{\R}$-CW complexes. By \cite[Thm B(a)]{Carlsson_Sullivans_Conj}, the natural map $\pi_0(\Alb_1^{G_{\R}}) \to \pi_0(\Alb_1^{h G_{\R}})$ is a bijection. Since $\Alb_1$ is a $\rK(\pi/[\pi]_2,1)$, the natural map $ \pi_0(\Alb_1^{h G_{\R}}) \to \mathscr{S}_{\pi_1(\Alb_1/G_{\R})}$ is a bijection. Under these bijections, $\alpha$ is identified with $\kappa^{\ab}$. The same reasoning applied to $\Alb_2$ identifies the image of $ \pi_0(\Alb_2^{G_{\R}}) \to \pi_0(\Alb_1^{G_{\R}})$ with the image of $\mathscr{S}_{\pi_1(\Alb_2/G_{\R})} \to \mathscr{S}_{\pi_1(\Alb_1/G_{\R})}$. This shows that Theorem \ref{Top_2nilSC} and Theorem \ref{2nilsectionconjecture} are equivalent.
\end{proof}

\epoint{Example}\label{Alb3example} Let $X=\pmR$ equipped with a real base point $b$ in $(0,1)$, say $b = \frac{1}{2}$. The fundamental group $\pi$ is freely generated by $x_1$ and $x_2$, where $x_1$ is represented by the loop $t \mapsto e^{2 \pi i t}/2$ for $t \in [0,1]$ and $x_2$ is the image of $x_1$ under the automorphism of $X$ given by $z \mapsto 1-z$. 

The set $\{ x_1, x_2, [x_1, x_2] \}$ is a basis for $\pi/[\pi]_2 \oplus [\pi]_2/[\pi]_3$, so points of $(\pi/[\pi]_2 \oplus [\pi]_2/[\pi]_3) \otimes_{\Z} \R$ can be labeled $(a_1, a_2, a_{12}) \in \R^3$ as an abbreviation for $a_1 x_1 + a_2 x_2 + a_{12} [x_1, x_2]$. The action of $\pi/[\pi]_3$ on $(\pi/[\pi]_2 \oplus [\pi]_2/[\pi]_3) \otimes_{\Z} \R$ is given by $x_1 (a_1, a_2, a_{12}) = (a_1+1, a_2, a_{12})$ and $x_2 (a_1, a_2, a_{12}) = (a_1, a_2 + 1, a_{12} -a_1)$. Note that $[x_1,x_2] (a_1, a_2, a_{12}) = (a_1, a_2, a_{12}+1)$ as well. It follows that $\Alb_3$ is the quotient of the unit cube in $\R^3$ given by identifying the $a_1 = 0$ face with the $a_1 = 1$ face via the translation $x_1$, identifying the $a_{12} = 0$ face with the $a_{12} = 1$ face via the translation $[x_1,x_2]$, and identifying the $a_2 = 0$ face with the $a_2 = 1$ face via the translation-shear $x_2$. The $G_{\R}$-action on $\Alb_3$ is given by $(a_1,a_2,a_{12}) \mapsto (-a_1, -a_2, a_{12})$. 

$\Alb_2$ is the torus given as the quotient of the unit square in $(\pi/[\pi]_2) \otimes_{\Z} \R$ with respect to the basis $\{x_1, x_2 \}$ by the translations $(a_1, a_2) \mapsto (a_1 + 1, a_2)$ and $(a_1, a_2) \mapsto (a_1, a_2  + 1)$, and $G_{\R}$ acts by multiplication by $-1$.  

$X$ deformation retracts $G_{\R}$-equivariantly onto the union of the two circles which are the images of the representative loops for $x_1$ and $x_2$ described above. $\alpha_2$ is the map taking the point $e^{2 \pi i t}/2$ in the image of $x_1$ to $t x_1$ in $\Alb_3$ and similarly for the points in the image of $x_2$. 

Note that there are four $G_{\R}$-fixed points of $\Alb_2$ given by the $2$-torsion points $$\{ (0,0), (\frac{1}{2}, 0), (0, \frac{1}{2}), (\frac{1}{2}, \frac{1}{2}) \}$$ of $\Alb_2$, and that the first three lift to fixed points of $\Alb_3$ as they constitute the image of the fixed points of $X$ in $\Alb_2$. The fourth point $(\frac{1}{2}, \frac{1}{2})$ does not lift to a fixed point of $\Alb_3$ as $G_{\R}$ acts on the fiber above $(\frac{1}{2}, \frac{1}{2})$ by translation by $(0,0,\frac{1}{2})$. The fact that $(\frac{1}{2}, \frac{1}{2})$ does not lift can also be seen by applying Theorem \ref{Top_2nilSC}.  Example \ref{Alb3example} is illustrated with figure \ref{figureAlb}. The fixed points are shown in red.

\begin{figure} [ht]
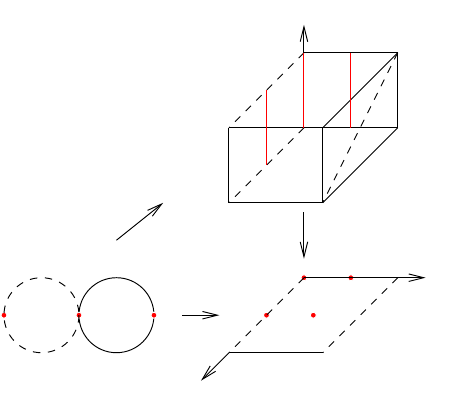
\caption{Approximations of $\pmR$}
\label{figureAlb}
\end{figure}

\appendix \section{Symmetric powers of smooth curves over their generalized Jacobians}\label{appendixa}

We provide a proof of Proposition \ref{SymnPicn} in the case where $X$ is non-proper. See \cite[VII \S 2 Prop 2.1]{Arbarello_I} for the case where $X$ is proper. This proposition and the following proof are both well-known, but we were unable to locate a reference in the literature.

For clarity, consider the more general situation in which $Y$ is a geometrically integral, proper curve over $k$ equipped with a $k$-point $y_0$.  Let $\Div^n_{Y}$ denote the functor taking a locally Noetherian scheme $T$ over $k$ to the closed subschemes $D$ of $Y \times T$, flat over $T$ of degree $n$, and with invertible ideal sheaf $\mathcal{I}$ (see \cite[8.2 p 212]{BLR}). The association $D \mapsto \mathcal{I}^{-1}$ for $D \in \Div_{Y}(T)$ determines a map \begin{equation}\label{YDivntoPicneq}\Div^n_{Y} \rightarrow \Pic^n_{Y}.\end{equation} Let $\Div^n_{Y,y_0}$ be the functor taking $T$ to the subset of $\Div^n_{Y}(T)$ consisting of those closed subschemes containing $y_0 \times T$. The map $\Div^n_{Y,y_0} \to \Div^n_{Y}$ is represented as follows. 

By \cite[8.1 Prop. 4]{BLR}, there is a unique universal invertible sheaf $\mathcal{P}$ over $Y \times \Pic_{Y}$ whose restriction to $y_0 \times \Pic_{Y}$ is trivial.  Let $p: Y \times \Pic_{Y} \rightarrow \Pic_{Y}$ and $p': Y \times \Pic_{Y} \rightarrow Y$ denote the projections. The closed immersion $y_0: \Spec k \to Y$ corresponds to a short exact sequence of sheaves on $Y$ $$0 \to \mathcal{I}_{y_0} \to \mathcal{O}_{Y} \to \mathcal{O}_{y_0} \to 0.$$ Applying $(p')^*$ and tensoring with $\mathcal{P}$ yields the short exact sequence \begin{equation}\label{Potimesp'*inftyseq}0 \rightarrow \mathcal{P} \otimes (p')^* \mathcal{I}_{y_0} \to \mathcal{P} \to \mathcal{P} \otimes (p')^*  \mathcal{O}_{y_0} \to 0. \end{equation} By \cite[8.1 Thm 7]{BLR}, since $\mathcal{P}$ is flat over $\Pic_{Y}$, there is a coherent sheaf $\mathcal{F}$ on $\Pic_{Y}$ with functorial isomorphisms $$p_*( \mathcal{P} \otimes p^* \mathcal{M}) \cong \sHom (\mathcal{F}, \mathcal{M})$$ for all quasi-coherent sheaves $\mathcal{M}$ on $\Pic_{Y}$. The morphism $\mathcal{P} \to \mathcal{P} \otimes (p')^*  \mathcal{O}_{y_0}$ defines a natural map \begin{equation}\label{PtoPy_0_yields}p_*( \mathcal{P} \otimes p^* \mathcal{M}) \to p_*(  \mathcal{P} \otimes (p')^*  \mathcal{O}_{y_0} \otimes p^* \mathcal{M}) \cong \mathcal{M},\end{equation} where the isomorphism comes from the trivialization of $\mathcal{P}$ restricted to $y_0 \times \Pic_Y$.  The map \eqref{PtoPy_0_yields} defines a map $\mathcal{O}_{\Pic_Y} \to \mathcal{F}.$ Let $\mathcal{E}$ denote the cokernel. 

For a coherent sheaf $\mathcal{M}$, let $\proj (\mathcal{M}) = \Sheafproj \Sym \mathcal{M}$ as in \cite[3.1.3]{egaII}. When $\mathcal{M}$ is locally free, $\proj (\mathcal{M})$ is called a projectivized vector bundle.

The proof of \cite[8.2 Prop 7]{BLR} can be modified to show:

\tpoint{Proposition}\label{DivYyo_to_DivY_rep}{\em $\Div_{Y, y_0} \hookrightarrow \Div_{Y}$ is represented by the closed immersion $\proj \mathcal{E} \to \proj \mathcal{F}$.}

\begin{proof}
It is sufficient to show the claim after restricting $\Div_{Y,y_0}$ and $\Div_{Y}$ to functors on schemes $T \rightarrow \Pic_{Y}$.  We will let $p$, $p'$ also denote their pullbacks to $Y \times T$, $\mathcal{P}_T$ denote the pullback of $\mathcal{P}$ to $Y \times T$, and $\mathcal{F}_T$, and $\mathcal{E}_T$ denote the pullbacks to $T$ of $\mathcal{F}$, and $\mathcal{E}$ respectively. For a point $t: \Spec k(t) \rightarrow T$, let $Y_{k(t)} = (Y \times T) \times_T \Spec k(t)$. For an invertible sheaf $\mathcal{L}$, let $\mathcal{L}^{-1}$ be the dual invertible sheaf $\mathcal{L}^{-1} = \Hom_{\mathcal{O}}(\mathcal{L}, \mathcal{O})$.

A section $s: \mathcal{O}_{Y \times T} \rightarrow \mathcal{P}_T$ induces a map $$s^{-1}: \mathcal{P}_T^{-1} \rightarrow \mathcal{O}_{Y \times T}^{-1} = \mathcal{O}_{Y \times T}$$ which is an injection such that the corresponding closed subscheme $D$ is flat over $T$ if and only if the restriction of $s^{-1}$ to $Y_{k(t)}$ is injective for all $t$ by \cite[Prop. 11.3.7]{egaIV_3}. Since $Y_{k(t)}$ is reduced and irreducible, the restriction of $s^{-1}$ to $Y_{k(t)}$ is injective if and only if it is non-zero.  By the definition of $\mathcal{F}$, the set of such $s$ is in natural bijection with morphisms $\mathcal{F}_T \rightarrow \mathcal{O}_T$ which are non-zero on all stalks. By Nakayama's lemma, a morphism $\mathcal{F}_T \rightarrow \mathcal{O}_T$ is non-zero on all stalks if and only if it is surjective.

$D$ contains $y_0 \times T$ if and only if the image of $s^{-1}$ is contained in the ideal sheaf $(p')^* \mathcal{I}_{y_0}$ of $y_0 \times T$. This occurs if and only if the image of $s$ is contained in $\mathcal{P}_T \otimes (p')^* \mathcal{I}_{y_0}$, which by \eqref{Potimesp'*inftyseq} occurs if and only if the composition of $s$ with $$\mathcal{P}_T \to \mathcal{P} \otimes (p')^*  \mathcal{O}_{y_0} $$ is $0$. Thus, $D$ contains $y_0 \times T$ if and only if the morphism $\mathcal{F}_T \rightarrow \mathcal{O}_T$ corresponding to $s$ is pulled back from a morphism $\mathcal{E}_T \rightarrow \mathcal{O}_T$.  

It follows that the association $s \mapsto D$ induces a bijection between the set of elements $D$ of $\Div_{Y, y_0}$ equipped with an isomorphism $\mathcal{I} \rightarrow \mathcal{P}_T^{-1}$, where $\mathcal{I}$ denotes the ideal sheaf of $D$, and surjections $\mathcal{E}_T \rightarrow \mathcal{O}_T$. Let $\Hom(\mathcal{E}_T, \mathcal{O}_T)^{\textrm{surj}} \subset \Hom(\mathcal{E}_T, \mathcal{O}_T)$ denote the subsheaf of surjections. As the automorphisms of $\mathcal{P}_T^{-1}$ are canonically isomorphic to the global sections $\Gamma (Y \times T, \mathcal{O}_{Y \times T}^*) = \Gamma(T, \mathcal{O}_T^*)$\hidden{$p_*(\mathcal{O}_{Y \times T}) = \mathcal{O}_T$ by \cite[7.8.6]{egaIII_2}}, the association $s \mapsto D$ induces a bijection between the set of elements $D$ of $\Div_{Y, y_0}$ such that $\mathcal{I} \cong \mathcal{P}_T^{-1}$ and $\Gamma(T,\Hom(\mathcal{E}_T, \mathcal{O}_T)^{\textrm{surj}})/\Gamma(T, \mathcal{O}^*)$.

Since two invertible sheaves $\mathcal{L}_1$ and $\mathcal{L}_2$ on $Y \times T$ induce the same map $T \rightarrow \Pic_{Y}$ if and only if $\mathcal{L}_1 \otimes \mathcal{L}_2^{-1}$ is pulled back from $T$ \cite[8.1 Prop 4]{BLR}, it follows that $$\Gamma(T, \Hom(\mathcal{E}_T, \mathcal{O}_T)^{\textrm{surj}}/\mathcal{O}_T^*) = (\Div_{Y} \times_{\Pic_{Y}} T )(T).$$ Note that $\Gamma(T, \Hom(\mathcal{E}_T, \mathcal{O}_T)^{\textrm{surj}}/\mathcal{O}_T^*)$ is in natural bijection with the set of all equivalence classes of pairs of an invertible sheaf $\mathcal{L}$ on $T$ and a surjection $\varphi: \mathcal{E}_T \rightarrow \mathcal{L}$ where $(\mathcal{L}, \varphi)$ and $(\mathcal{L}', \varphi' )$ are equivalent if there is an isomorphism $\theta: \mathcal{L} \rightarrow \mathcal{L}'$ such that $\varphi' = \theta \varphi$.  By \cite[4.2.3]{egaII}, it follows that the canonical closed immersion  $\proj (\mathcal{E}) \rightarrow \proj (\mathcal{F})$ represents $\Div_{Y, y_0} \rightarrow \Div_{Y}$. 
\end{proof}

Let $f$ be the global section of $\mathcal{F}$ corresponding to the map $\mathcal{O}_{\Pic_Y} \to \mathcal{F}$. By Proposition \ref{DivYyo_to_DivY_rep}, \begin{equation}\label{Div-Div0-sheafspec}\Div_{Y} - \Div_{Y, y_0} \cong \Sheafspec_{\Pic Y} \Sym \mathcal{F} [f^{-1}]_0,\end{equation} where $\Sym \mathcal{F} [f^{-1}]_0$ denotes the degree $0$ elements of the graded algebra given by inverting $f$ in the symmetric algebra of $\mathcal{F}$. We will show that for large enough $n$, $(\Div^n_{Y} - \Div^n_{Y, y_0}) \to \Pic^n Y$ is an affine bundle, i.e. $\Sym \mathcal{F} [f^{-1}]_0$ is locally isomorphic to a polynomial algebra over $\mathcal{O}_{\Pic^n_Y}$ with affine transition maps between the local isomorphisms. For this, we will need the following lemma.

\tpoint{Lemma}\label{cohsheavotimesPcohflat}{\em For a coherent sheaf $\mathcal{M}$ on $Y$, let $\mathcal{G}$ be the coherent sheaf $(p')^* \mathcal{M} \otimes \mathcal{P}$ on $Y \times \Pic_{Y}$. Then: 
\begin{enumerate}% \item\label{Gflat} $\mathcal{G}$ is flat over $\Pic_{Y}$.  
\item \label{Rbig0} For large $n$, the restriction of $R^i p_* \mathcal{G}$ to $\Pic^n_{Y}$ is $0$ for all $i>0$, and locally free for $i=0$.
\item \label{Gcohflat0} For large $n$, the restriction of $\mathcal{G}$ to $ Y \times \Pic_{Y}^n$ is cohomologically flat over $\Pic^n_{Y}$ in all dimensions.
\end{enumerate}}

\begin{proof}
%\eqref{Gflat}:  Near any point of $Y \times \Pic_{Y}$, note that $\mathcal{G}$ is isomorphic to $(p')^* \mathcal{M}$, which is flat because $\mathcal{M}$ is flat over $\Spec k$.

The $i>0$ case of \eqref{Rbig0}: Since $Y$ is projective, there is a relatively very ample invertible sheaf $\mathcal{L}$ for $p:Y \times \Pic_{Y} \to \Pic_{Y}$. Let $ \mathcal{L}^{\otimes n}$ denote the $n$-fold tensor product of $\mathcal{L}$. Let $\mathcal{G}^k$ denote the restriction of $\mathcal{G}$ to $Y \times \Pic_{Y}^k$. For a fixed $k$, there is $N$ such that for $n >N $, we have $R^ip_* (\mathcal{G}^k \otimes  \mathcal{L}^{\otimes n}) = 0$ for all $i>0$ by \cite[Thm 2.2.1]{egaIII_1}. The invertible sheaf $\mathcal{P} \otimes \mathcal{L}^{\otimes m}$ induces an isomorphism $t: \Pic_{Y}^k \rightarrow \Pic_{Y}^{k+md}$ such that $t^* \mathcal{G}^{k+md}= \mathcal{G}^k \otimes \mathcal{L}^{\otimes m}$, where $d$ denotes the degree of $\mathcal{L}$. Since  $t^* R^ip_* (\mathcal{G}^{k+md}) = R^i p_* (t^* \mathcal{G}^{k+md})$, we have that $R^ip_* \mathcal{G}^{k+md} = 0$ for a fixed $k$, and all $m$ sufficiently large and $i>0$. Taking $k=0,1,\ldots ,d-1$ shows that for $m$ sufficiently large we have $R^ip_* \mathcal{G}^{m} = 0$ for all $i>0$.

\eqref{Gcohflat0}: By the $i>0$ case of \eqref{Rbig0}, the restriction to $Y \times \Pic^n_{Y}$ of $R^ip_* (\mathcal{G})$ for $n$ sufficiently large is locally free for all $i \geq 1$. By \cite[Prop 7.8.5]{egaIII_2} it follows that this restriction of $\mathcal{G}$ is cohomologically flat in dimensions $i \geq 1$. For a point $z$ of $\Pic_{Y}$, let $\mathcal{G}_z$ denote the pullback of $\mathcal{G}$ by the closed immersion $Y_{k(z)}= Y \times \Spec k(z) \rightarrow Y \times \Pic^n_{Y}$ corresponding to $z$. By \cite[Prop 7.8.4]{egaIII_2}, it follows that for a fixed  for $i\geq 1$, the function $z \mapsto d_i (z) = \operatorname{dim}_{k(z)} \rH^i(Y_{k(z)}, \mathcal{G}_z)$ on the points of $\Pic^n_{Y}$ for $n$ sufficiently large is locally constant. Since the Euler characteristic of $\mathcal{G}_z$ is locally constant \cite[Thm 7.9.4]{egaIII_2}, it follows that $d_0$ is locally constant when restricted to $\Pic^n_{Y}$ for $n$ sufficiently large. Since $\Pic^n_{Y}$ is reduced,  \eqref{Gcohflat0} follows from \cite[Prop 7.8.4]{egaIII_2}. 

The $i=0$ case of \eqref{Rbig0}, follows from \eqref{Gcohflat0} and \cite[7.8.5]{egaIII_2}. 
\end{proof}

By \cite[8.1 Thm 7]{BLR}, when $n$ is large enough so that $\mathcal{P}$ is cohomologically flat over $\Pic_Y^n$ in dimension $0$, the dual of $p_*( \mathcal{P})$ is canonically isomorphic to $\mathcal{F}$. Restrict \eqref{Potimesp'*inftyseq} to $Y \times \Pic^n Y$ and apply $p_*$. By Lemma \ref{cohsheavotimesPcohflat}, we obtain the short exact sequence of locally free sheaves \begin{equation}\label{p_*Potimes_etc} 0 \rightarrow p_*( \mathcal{P} \otimes (p')^* \mathcal{I}_{y_0}) \to p_*( \mathcal{P}) \to \mathcal{O}_{\Pic_Y^n} \to 0,\end{equation} for sufficiently large $n$.  Dualizing yields the short exact sequence of locally free sheaves $$0 \to \mathcal{O}_{\Pic_Y^n} \to \mathcal{F} \to \mathcal{E} \to 0.$$ By \eqref{Div-Div0-sheafspec}, it follows that for large enough $n$, $(\Div^n_{Y} - \Div^n_{Y, y_0}) \to \Pic^n Y$ is an affine bundle.

By \cite[6.3.9]{sga4III}, $\Sym^n X$ represents the functor taking a locally Noetherian scheme $T$ over $k$ to the closed subschemes $D$ of $X \times T$ which are flat and finite over $T$ of degree $n$ with invertible ideal sheaf. Thus $\Sym^n X \cong \Div^n_{X^+} - \Div^n_{X^+, \infty}$, where $\infty$ is the unique $k$ point in $X^+ - X$. It follows that $\Sym^n X \to \Pic^n X^+$ is an affine bundle, which also shows Proposition \ref{SymnPicn}. 

} % end of parskip; it started just before the introduction

\bibliographystyle{delta2real2}

%the bibliographystyle was modified by Joe from amsalpha. If I change bibliographystyle,  "thesis" bibliography will still compile, but without the overides to cite EGA.

\bibliography{delta2real2}

%LaTeX BibTeX LaTeX  LaTeX 
%\cite[pg 27] {Ihara_GT}

\end{document}